\documentclass{article}

\usepackage{bbm}
\usepackage{amsmath}
\usepackage{amsfonts}
\usepackage{amssymb}
\usepackage[mathscr]{eucal}
\numberwithin{equation}{section}
\usepackage{amsthm}
\newtheorem{theorem}{Theorem}

\newtheorem{remark}[theorem]{Remark}

\usepackage[a4paper,left=3.18cm, right=3.18cm, top=2.54cm, bottom=2.54cm]{geometry}
\usepackage{enumitem}
\setlist[enumerate]{leftmargin=.5in}
\setlist[itemize]{leftmargin=.5in}

\usepackage{graphicx} 
\graphicspath{{./}}
\usepackage{multirow}
\usepackage{booktabs}
\usepackage{subcaption}
\newlength{\Oldarrayrulewidth}
\newcommand{\Cline}[2]{
  \noalign{\global\setlength{\Oldarrayrulewidth}{\arrayrulewidth}}
  \noalign{\global\setlength{\arrayrulewidth}{#1}}\cline{#2}
  \noalign{\global\setlength{\arrayrulewidth}{\Oldarrayrulewidth}}}
\usepackage{algorithm}
\usepackage{algpseudocode}
\usepackage{xcolor}
\definecolor{trueblue}{rgb}{0.0, 0.45, 0.81}
\usepackage[colorlinks=true, allcolors=trueblue]{hyperref}

\usepackage{titlesec}
\titleformat*{\section}{\Large\bfseries\sffamily}
\titleformat*{\subsection}{\large\bfseries\sffamily}

\renewenvironment{abstract}{
  \quotation
  \textbf{\textsf{{\abstractname.}}}
  \sffamily
}{\endquotation}

\usepackage{authblk}
\newcommand*\samethanks[1][\value{footnote}]{\footnotemark[#1]}
\makeatletter
\def\@fnsymbol#1{\ensuremath{\ifcase#1\or \dagger\or \ddagger\or
   \mathsection\or \mathparagraph\or \|\or **\or \dagger\dagger
   \or \ddagger\ddagger \else\@ctrerr\fi}}
\makeatletter

\def\0{\mathbf{0}}
\def\1{\mathbf{1}}

\def\c{\mathbf{c}}
\def\f{\mathbf{f}}
\def\bf{\mathbbm{f}}
\def\g{\mathbf{g}}

\def\h{\mathbf{h}}
\def\bh{\mathbbm{h}}
\def\n{\mathbf{n}}

\def\p{\mathbf{p}}
\def\r{\mathbf{r}}

\def\x{\mathbf{x}}
\def\y{\mathbf{y}}
\def\z{\mathbf{z}}

\def\rS{\mathrm{S}}
\def\rA{\mathrm{A}}

\def\R{\mathbb{R}}
\def\A{\mathcal{A}}

\def\bbA{\mathbb{A}}
\def\E{\mathcal{E}}
\def\F{\mathcal{F}}
\def\bF{\mathtt{F}}

\def\M{\mathcal{M}}

\def\bN{\mathtt{N}}

\def\T{\mathcal{T}}
\def\V{\mathcal{V}}
\def\bV{\mathscr{V}}

\def\bC{\mathbf{C}}

\def\bY{\mathbf{Y}}
\def\d{\mathrm{d}}

\def\bphi{{\boldsymbol{\phi}}}

\def\btheta{{\boldsymbol{\theta}}}

\newcommand{\argmin}{\operatornamewithlimits{argmin}}

\newcommand{\SD}{\operatornamewithlimits{SD}}

\title{ 
\bfseries\sffamily{
Spherical Area-Preserving Parameterization via Energy Minimization}

}
\author{
\sffamily
{Shu-Yung Liu\thanks{\footnotesize Department of Mathematics, National Taiwan Normal University, Taipei, Taiwan (\href{mailto:lii227857@gmail.com}{lii227857@gmail.com}, \href{mailto:yue@ntnu.edu.tw}{yue@ntnu.edu.tw}) }
\, and \, Mei-Heng Yueh\samethanks[1] }  
}

\date{}

\begin{document}
\captionsetup[figure]{labelfont={bf,sf},name={Figure},labelsep=period}
\captionsetup[table]{labelfont={bf,sf},name={Table},labelsep=period}
\maketitle

\begin{abstract}
We propose a novel method, called {\em spherical authalic energy minimization} (SAEM), for computing spherical area-preserving parameterizations of genus-zero closed surfaces, with strong theoretical foundations. The global convergence of the associated computational algorithm is theoretically guaranteed. In addition, we introduce a Riemannian bijective correction method that ensures the bijectivity of the resulting mapping under mild assumptions. Numerical experiments show that SAEM effectively minimizes area distortion and achieves bijective mappings, outperforming state-of-the-art methods. Finally, we demonstrate the practical utility of SAEM in shape description.


\end{abstract}

\section{Introduction}
Surface parameterization involves mapping a surface in 3-dimensional space onto a simpler canonical domain (e.g., a disk or sphere). This technique has been used extensively in computer graphics for tasks such as surface resampling and registration \cite{LaLu14,LuLY14,YoMY14,YuLW17}, and it has also demonstrated its utility in biomedical applications, including brain flattening and tumor segmentation \cite{AnSH99, AnHT99,LiJY21,LiLH22}. For a more comprehensive overview and classical applications, readers may refer to the survey papers \cite{FlHo05,ShPR06} and the lecture notes \cite{HoLe07}.

A desirable parameterization preserves the geometric features of a surface as much as possible. In practice, achieving an isometric (length-preserving) mapping is generally infeasible. Consequently, most approaches focus on either angle preservation (conformal mappings) or area preservation (authalic mappings), while some attempt to balance both in pursuit of an approximate isometry.

For simply connected closed surfaces, the unit sphere is a common parameter domain, motivating numerous studies on spherical conformal parameterizations. Haker et al. \cite{AnSH99, AnHT99, HaAT00} used the finite element approach to compute spherical conformal mappings by solving partial differential equations on the extended complex plane. However, their approach, which depends on stereographic projection, results in significant angular distortion near the north pole. To remedy this issue, Choi et al. \cite{ChLL15} fixed the southern region and applied a specific quasi-conformal mapping to offset the distortion in the northern region. Likewise, Choi et al. \cite{ChHo16} and Yueh et al. \cite{YuLL19, LiHL23} addressed angular distortion in the northern region by applying the south-pole stereographic projection and iterating the process until convergence. Alternatively, Nadeem et al. \cite{NaSZ17} applied Haker et al.’s method separately to the northern and southern hemispheres, then conformally welded the two solutions to minimize angular distortion.

While angle-preserving mappings maintain the local shape of a surface, they often result in significant area distortion. Area preservation is crucial for tasks such as brain morphometry. Several studies have focused on spherical area-preserving mappings for simply connected closed surfaces. Moser \cite{Mose65} introduced a method by solving two differential equations in the sphere. Building on this, Dominitz et al. \cite{DoTa10} computed optimal mass transportation (OMT) mappings through a gradient-flow method. Subsequent studies extensively applied the OMT theory. Nadeem et al. \cite{NaSZ17} computed OMT mappings by minimizing a convex energy in Euclidean space via stereographic projection. Similarly, Choi et al. \cite{ChGK22} computed OMT mappings in Euclidean space but first punctured a quadrilateral region at the bottom of the surface. In contrast, Cui et al. \cite{CuQW19} developed a computational algorithm to compute OMT mappings directly on the sphere. 

In addition to the OMT approaches, Lyu et al. \cite{LyCh24} employed density-equalizing methods combined with tangent plane projection and retraction to achieve area preservation. Meanwhile, Yueh et al. \cite{YuLL19} and Huang et al. \cite{HuLW24} introduced stretch energy minimization (SEM), using a combination of stereographic projection and the fixed-point method \cite{YuLW19}. The theoretical link between stretch energy and area preservation was established in \cite{Yueh23}. This approach was further enhanced by incorporating a Riemannian gradient descent method with theoretical convergence \cite{SuYu24}.

From an application perspective, spherical authalic parameterization is essential for shape description with spherical harmonics. Spherical harmonics \cite{Byer59} are special functions defined on the unit sphere, forming an orthonormal basis. They are widely used as shape descriptors \cite{BrGK95} to effectively capture surface information with appropriate coefficients. This results in interest in medical applications, including analyzing brain ventricle morphology in twins \cite{GeSJ01}, magnetic resonance imaging in schizophrenia \cite{StLP04}, and hippocampal shape features in Alzheimer's disease \cite{GeCC09, ChLS20}.

Recently, for simply connected open surfaces, Liu and Yueh \cite{LiYu24} refined the energy minimization approach, significantly improving area preservation compared to the OMT method \cite{ZhSG13} and density-equalizing map \cite{ChRy18}. However, this method has not been well developed for genus-zero closed surfaces. In this paper, we adapt and modify this energy minimization approach for spherical area-preserving parameterization. We expect it will outperform current state-of-the-art methods, as in the open-surface setting, and be suitable for shape description via spherical harmonics.

The paper is organized as follows. In Section \ref{sec:1.2}, we highlight the contributions of our work. Section \ref{sec:2} provides the mathematical background. In Section \ref{sec:3}, we introduce the modified objective functional for improved numerical performance, with the associated iterative method detailed in Section \ref{sec:4}. Section \ref{sec:5} discusses a post-processing step to ensure bijectivity for resulting mappings. Numerical results and applications in shape description via spherical harmonics are presented in Sections \ref{sec:6} and \ref{sec:7}, respectively. Finally, we discuss our proposed method compared with existing methods in Section \ref{sec:8} and conclude in Section \ref{sec:9}.

\section{Contributions}
\label{sec:1.2}
The main contributions of this work are as follows:
\begin{itemize}
\item[(i)] \textbf{Novel Spherical Authalic Energy}: 
We propose a novel objective functional, the spherical authalic energy, which improves bijectivity during optimization for spherical area-preserving parameterization of genus-zero closed surfaces.
\item[(ii)] \textbf{Efficient Energy Minimization with Guaranteed Convergence}: 
We develop an energy minimization method with theoretically guaranteed convergence, outperforming state-of-the-art methods in both accuracy and computational efficiency.
\item[(iii)] \textbf{Riemannian Bijective Correction Method}: 
We introduce a Riemannian bijective correction method with strong theoretical support, which is capable of resolving hundreds of folded triangles, effectively addressing challenging cases.
\item[(iv)] \textbf{Practical Applications in Shape Description}:
We demonstrate the practical application of our method in shape description to show its utility and effectiveness.
\end{itemize}

\section{Mathematical background}
\label{sec:2}
In this section, we introduce some basic mathematical concepts closely related to our work. 

\subsection{Simplicial mapping}
A smooth surface can be approximated by a simplicial surface $\M$, A \textit{triangular mesh} is defined as 
$$
\T(\M) = \big(\V(\M), \E(\M), \F(\M) \big),
$$
where
\begin{align*}
\V(\M) &= \left\{ v_\ell = (v_\ell^1, v_\ell^2, v_\ell^3) \in\R^3 \right\}_{\ell=1}^n, \\[4pt]
\F(\M) &= \left\{ \tau_s = [ {v}_{i_s}, {v}_{j_s}, {v}_{k_s} ] \subset\R^3 \mid {v}_{i_s}, {v}_{j_s}, {v}_{k_s}\in\V(\M) \right\}_{s=1}^m, \\[4pt]
\E(\M) &= \left\{ [v_i,v_j] \subset\R^3 \mid [v_i,v_j,v_k]\in\F(\M)  \right\}.
\end{align*}
Here, $\V(\M)$, $\F(\M)$, and $\E(\M)$ are the sets of $n$ vertices, $m$ oriented triangular faces, and undirected edges, respectively. The bracket $[{v}_{i}, {v}_{j}, {v}_{k}]$ denotes the oriented \textit{$2$-simplex}, i.e., the triangle with those vertices.
The simplicial surface $\M$ is the \textit{geometric realization} of $\T(\M)$, which is the union of flat triangles in $\R^3$. i.e., $\M = \cup_{\tau\in\F(\M)} \tau$.

Each triangle is determined by three vertex positions associated with a face in $\F(\M)$.

Let $\mathbb{S}^2$ be the unit sphere in $\mathbb{R}^3$. A {\it spherical simplicial mapping} is a piecewise affine mapping $f: \M \to \mathbb{S}^2$ that sends every vertex of the simplicial surface $\M$ to a point on $\mathbb{S}^2$. 
Because $f$ is affine on each triangular face, the image $f(\M)$ 

itself is a simplicial surface lying on $\mathbb{S}^2$.

\begin{figure}
\centering
\includegraphics{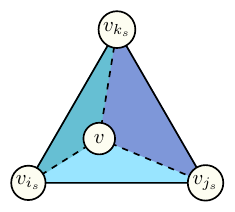}
\caption{An illustration of the barycentric coordinates on a triangular face.}
\label{fig:barycentric}
\end{figure}

The simplicial mapping is completely determined by its values on the vertices. Any point $v$ inside a triangular face $\tau = [{v}_{i}, {v}_{j}, {v}_{k}]$ can be written in barycentric form
$$
v = \lambda_i v_i + \lambda_j v_j + \lambda_k v_k, 
$$
with
$$
\lambda_i = \frac{|[v, v_j, v_k]|}{|\tau|}, \quad
\lambda_j = \frac{|[v_i, v, v_k]|}{|\tau|}, \quad
\lambda_k = \frac{|[v_i, v_j, v]|}{|\tau|},
$$
where {$|\tau|=\int_\tau \,\d A$} denotes the area of the triangle $\tau$ (see Figure~\ref{fig:barycentric}).

Then, the piecewise affine mapping $f$ is given by
$$
f_i \equiv f(v_i) = (f_i^1, f_i^2, f_i^3)^\top, \quad \text{for every $v_i\in\V(\M)$},
$$
and, on each face,
$$
f|_{\tau}(v) =  \lambda_i \, f(v_{i}) + \lambda_j \, f(v_{j}) + \lambda_k \, f(v_{k}),
$$
for every $\tau\in\F(\M)$.
{ 
A simplicial map $f$ is said to be authalic or area-preserving up to a global scaling if
\begin{equation} \label{eq:area-preserving}
|f(\tau)| 
= c \, |\tau|, ~~ \mbox{for all}~ \tau \in \mathcal{F(M)},
\end{equation}
where $c>0$ is a global constant. 
Note that the area of the image of $\tau$ under $f$ can be expressed as
$$
|f(\tau)| = \int_\tau \sqrt{\det( \mathrm{I}_f )} \,\d A,
$$
where $\mathrm{I}_f = {J}_f^\top {J}_f$ is the \textit{first fundamental form} of the map $f$, and $J_f = \begin{bmatrix}
\frac{\partial f}{\partial u} & \frac{\partial f}{\partial v} \end{bmatrix}$ denotes its Jacobian matrix in local coordinates. 
Therefore, the condition in \eqref{eq:area-preserving} is equivalent to requiring the area element $\sqrt{\det( \mathrm{I}_f )} = c$ to be constant across the surface. 
This property aligns with the classical concept of area-preserving maps in differential geometry.
}

\subsection{Authalic energy functional} \label{sec:2.2}
The \textit{authalic energy} is a global measure of how well a simplicial map $f$ preserves area. 

Introduced by Liu and Yueh \cite{LiYu24}, it is defined as
\begin{equation}
E_\rA(f) = \frac{|\M|}{\A(f)} E_\rS(f) - \A(f),
\label{eq:Ea}
\end{equation}
where { 
$$
|\M| = \int_\M \d A = \sum_{\tau\in\F(\M)} \int_\tau \d A = \sum_{\tau\in\F(\M)} |\tau|,
$$
is the area of $\M$, and
$$
\A(f) = \int_\M \sqrt{\det( \mathrm{I}_f )} \, \d A = \sum_{\tau\in\F(\M)} \int_\tau \sqrt{\det( \mathrm{I}_f )} \, \d A = \sum_{\tau\in\F(\M)} |f(\tau)|,
$$
is the total area of the mapped surface.
}
Here, $E_\rS$ denotes the \textit{stretch energy} defined by
\begin{equation}
E_\rS(f) = \sum_{\tau \in\mathcal{F(M)}} \frac{|f(\tau)|^2}{|\tau|}.
\label{eq:Es}
\end{equation}
The stretch energy was originally introduced in \cite{YuLW19} and later refined to the form~\eqref{eq:Es} in \cite{Yueh23}.

Under the assumption that $\A(f) = \bigl|\M\bigr|$, the lower bound of $E_\rS(f)$ is the image area $\A(f)$, and this minimum is attained when $f$ is area-preserving \cite[Theorem 3.3]{Yueh23}. By relaxing this assumption, the authalic energy $E_\rA(f)$ is introduced, which attains its minimum value of $0$ when $f$ is area-preserving \cite[Theorem 1]{LiYu24}. This relaxation allows the image area $\A(f)$ to vary during the optimization process.

\begin{figure}[t]
\centering
\includegraphics[height=5.5cm]{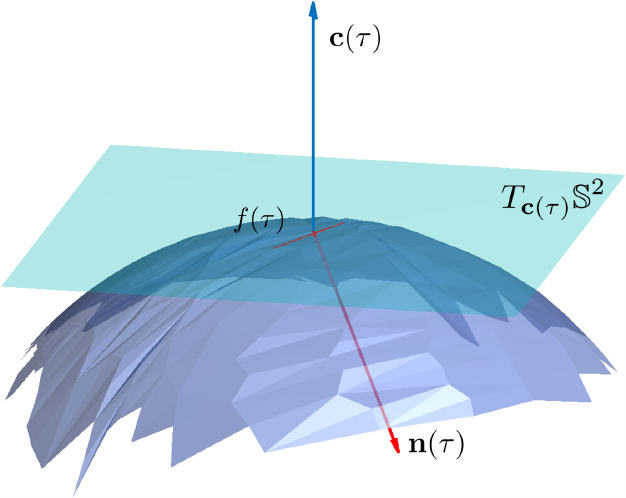}
\caption{ { 
Illustration of a spherical folded face $f(\tau)$ (red triangle), identified by the negative dot product between the normal vector $\n(\tau)$ and the face center $\c(\tau)$.
} }
\label{fig:SphereFolding}
\end{figure}

\section{Objective energy functional for spherical mappings}
\label{sec:3}
Authalic energy minimization has been shown to be efficient for computing area-preserving mappings in unit-disk parameterization for simply connected open surfaces \cite{LiYu24}. However, extending this approach to spherical parameterizations of genus-zero closed surfaces introduces the additional challenge of ensuring bijectivity. In this section, we discuss the bijective issue and present our modified objective functional, called the \textit{spherical authalic energy}.

\begin{figure}[t]
\centering
\includegraphics{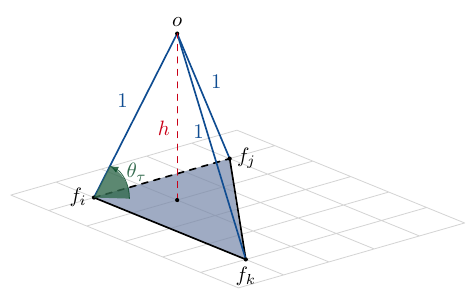}
\caption{The illustration of the tetrahedron by $[o, f_i, f_j, f_k]$ constructed by the triangle $f(\tau) = [f_i, f_j, f_k]$ and the origin $o$.}
\label{fig:tetrahedron}
\end{figure}

\subsection{Spherical authalic energy functional}
\label{sec:3.1}
{ 
Given a triangular face $\tau = [v_i, v_j, v_k]$, its face center is defined as
\begin{equation}
\c(\tau) = \frac{\widehat \c(\tau)}{\| \widehat \c(\tau) \|_2},~~~\widehat \c(\tau) = \frac{f(v_i) + f(v_j) + f(v_k)}{3},
\label{eq:FaceCenter}
\end{equation}
which corresponds to a tangent plane $T_{\c(\tau)}\mathbb{S}^2$. 
It's natural to define the {\em orientation of the face} $f(\tau)$ as the orientation of the projection of $f(\tau)$ in $T_{\c(\tau)}\mathbb{S}^2$, which can be identified as the sign of the dot product between $\c(\tau)$ and the normal vector of $f(\tau)$, denoted as $\n(\tau)$, i.e. $\mathrm{sign}(\n(\tau) \cdot \c(\tau))$. A face is {\em folded} if its orientation is negative.
Figure~\ref{fig:SphereFolding} illustrates a folded face $f(\tau)$, which is colored in red, whose normal vector $\n(\tau)$ is shown in red arrow. The blue arrow represents the face center $\c(\tau)$, and the corresponding tangent plane $T_{\c(\tau)}\mathbb{S}^2$ is shown in cyan. 

When we minimize the energy $E_\rA$ in \eqref{eq:Ea}, which contains a negative image area term $-\A(f)$ and an inverse dependence on $\A(f)$, the optimization tends to favor larger values of $\A(f)$ to reduce the overall energy. However, this often leads to the formation of folded triangles on the spherical mesh. 
In particular, when $\A(f)\geq 4\pi$, the image surface necessarily bounds a non-convex region in $\R^3$, which indicates the presence of orientation-reversing (i.e., folded) triangles. 
To correctly evaluate $\A(f)$ under such circumstances, it is necessary to compute the oriented area of each face, i.e., the signed face area.

Nevertheless, incorporating orientation into the definition of $\A(f)$ introduces discontinuities: when a triangle flips its orientation, the signed area of the face abruptly changes from positive to negative without passing through zero. Such discontinuities can significantly impact the optimization process.

}

To address this issue, we modify the energy functional as follows. Instead of using the face area of the spherical image $f(\tau) = [f_i, f_j, f_k]$, we consider the corresponding tetrahedron $[o, f_i, f_j, f_k]$, where $o$ is the origin, as illustrated in Figure~\ref{fig:tetrahedron}. As the orientation of the triangle $[f_i, f_j, f_k]$ on the sphere flips, the signed volume of $[o, f_i, f_j, f_k]$ continuously varies from a positive value to a negative one. This ensures that triangle orientations are accurately captured and that the functional remains continuous. 

Based on this concept, we define the {\em spherical authalic energy} functional as
\begin{equation} 
E_\bbA(f) = \frac{|\M|}{3\, \bV(f)} E_\rS(f) - 3\,\bV(f),
\label{eq:spherical Ea}
\end{equation}
where $E_\rS$ is the stretch energy defined in \eqref{eq:Es} and $\bV$ is a volume measure given by
\begin{equation} 
\bV(f) = \frac{1}{6} \sum_{[v_i, v_j, v_k] \in \mathcal{F(M)}} f_i^\top (f_j \times f_k).
\label{eq:volume measurement}
\end{equation}
Noting that the area of the unit sphere is $4\pi$, which is equal to three times its volume, the spherical authalic energy $E_\bbA$ in \eqref{eq:spherical Ea} can be viewed as an approximation to the authalic energy $E_\rA$ in \eqref{eq:Ea}. The following theorem provides a detailed statement.

\begin{theorem} \label{thm:sphere_Ea}
Let $\M$ be a closed genus-zero simplicial surface, and let $f:\M\to\mathbb{S}^{2}$ be an orientation-preserving simplicial map. Denote by $\varepsilon_f$ the maximum edge length of $f(\M)$. Assume that the image triangulation $f(\M)$ is shape-regular: there exists $\theta_0>0$ such that every interior angle of every triangle in $f(\M)$ is at least $\theta_0$. Set
\[
C=\frac{1}{4\sin^2\theta_0} \geq \frac{1}{3}.
\]
and assume that $C\varepsilon_f^2<1$. Then the spherical authalic energy \eqref{eq:spherical Ea} satisfies
\begin{equation}
\big|E_\bbA(f)-E_\rA(f)\big|
\leq \bigg(\A(f)+\frac{E_\rS(f)\,|\M|}{3\,\bV(f)}\bigg)
\Big(1-\sqrt{1-C\varepsilon_f^2}\Big),
\label{eq:4.3}
\end{equation}
where $E_\rA$ and $E_\rS$ are the authalic and stretch energies in \eqref{eq:Ea} and \eqref{eq:Es}, respectively, and $\bV$ is defined by \eqref{eq:volume measurement}. Moreover, if $f$ is area-preserving, then
\begin{equation}
0\leq E_\bbA(f)
\leq \bigg(\A(f)+\frac{\A(f)^2}{3\,\bV(f)}\bigg)
\Big(1-\sqrt{1-C\varepsilon_f^2}\Big).
\label{eq:4.4}
\end{equation}
\end{theorem}

\begin{proof}
By the definitions of $E_\rA$ and $E_\bbA$,
\begin{align}
E_\bbA(f)-E_\rA(f)
&=\bigg(\frac{|\M|}{3\,\bV(f)}E_\rS(f)-3\,\bV(f)\bigg)
 -\bigg(\frac{|\M|}{\A(f)}E_\rS(f)-\A(f)\bigg) \nonumber\\
&=E_\rS(f)|\M|\bigg(\frac{1}{3\,\bV(f)}-\frac{1}{\A(f)}\bigg)
 +\Big(\A(f)-3\,\bV(f)\Big) \nonumber\\
&=\bigg(1+\frac{E_\rS(f)\,|\M|}{3\,\bV(f)\,\A(f)}\bigg)
  \Big(\A(f)-3\,\bV(f)\Big).
\label{eq:4.5}
\end{align}

Fix a face $\tau\in\F(\M)$ with $f(\tau)=[f_i,f_j,f_k]$, and let the corresponding signed height $h_\tau= f_i^\top \n_\tau$. Since $f$ is orientation-preserving, $h_\tau\geq0$, and the corresponding tetrahedron satisfies
\begin{equation}
\frac{1}{6}\det[f_i,f_j,f_k] =\frac{h_\tau}{3}|f(\tau)|.
\label{eq:4.6}
\end{equation}
Let $p_\tau$ be the foot of the perpendicular from $o$ to this plane. Because $f_i,f_j,f_k\in\mathbb{S}^2$, Pythagoras' theorem gives
\[
|p_\tau-f_i|=|p_\tau-f_j|=|p_\tau-f_k|=\sqrt{1-h_\tau^2}.
\]
Thus $p_\tau$ is the circumcenter of $f(\tau)$, and its circumradius is
\[
R_\tau=\sqrt{1-h_\tau^2}.
\]

For an edge $[f_i,f_j]$ of $f(\tau)$, let $\theta_{i,j}^k(f)$ denote the interior angle at $f_k$ opposite that edge. Since the three interior angles sum to $\pi$ and each is at least $\theta_0$, we have $\theta_0\leq\pi/3$ and $\theta_{i,j}^k(f)\in[\theta_0,\pi-2\theta_0]$. Hence, $\sin\theta_{i,j}^k(f)\geq\sin\theta_0$. The circumradius formula therefore yields
\[
R_\tau
=\frac{|f_i-f_j|}{2\sin\theta_{i,j}^k(f)}
\leq\frac{\varepsilon_f}{2\sin\theta_0}.
\]
Consequently, $R_\tau^2\leq C\varepsilon_f^2$ and $C \geq \frac{1}{3}$, thus,
\begin{equation}
h_\tau=\sqrt{1-R_\tau^2}
\geq\sqrt{1-C\varepsilon_f^2}.
\label{eq:4.7}
\end{equation}
Using \eqref{eq:4.6} and \eqref{eq:4.7}, we obtain
\begin{align}
\A(f)-3\,\bV(f)
&=\sum_{\tau\in\F(\M)}|f(\tau)|(1-h_\tau) \nonumber\\
&\leq\A(f)\Big(1-\sqrt{1-C\varepsilon_f^2}\Big).
\label{eq:4.8}
\end{align}
Moreover, $0\leq h_\tau\leq1$ for every $\tau$, so $0\leq3\,\bV(f)\leq\A(f)$. It follows from \eqref{eq:4.5} that
\[
E_\bbA(f)-E_\rA(f)\geq0.
\]
Combining this observation with \eqref{eq:4.5} and \eqref{eq:4.8} gives
\begin{align*}
\big|E_\bbA(f)-E_\rA(f)\big|
&=\bigg(1+\frac{E_\rS(f)\,|\M|}{3\,\bV(f)\,\A(f)}\bigg)
  \Big(\A(f)-3\,\bV(f)\Big)\\
&\leq\bigg(\A(f)+\frac{E_\rS(f)\,|\M|}{3\,\bV(f)}\bigg)
  \Big(1-\sqrt{1-C\varepsilon_f^2}\Big),
\end{align*}
which proves \eqref{eq:4.3}.

If $f$ is area-preserving, then $E_\rA(f)=0$ and
$E_\rS(f)=\A(f)^2/|\M|$. Substituting these identities into \eqref{eq:4.3} yields
\begin{align*}
0\leq E_\bbA(f)
&\leq\bigg(\A(f)+\frac{E_\rS(f)\,|\M|}{3\,\bV(f)}\bigg)
  \Big(1-\sqrt{1-C\varepsilon_f^2}\Big)\\
&=\bigg(\A(f)+\frac{\A(f)^2}{3\,\bV(f)}\bigg)
  \Big(1-\sqrt{1-C\varepsilon_f^2}\Big),
\end{align*}
which proves \eqref{eq:4.4}.
\end{proof}

It is important to note that although $E_\bbA$ approximates $E_\rA$, the convergence behavior of minimizing $E_\bbA$ is significantly different from that of $E_\rA$. The numerical comparison of mappings produced by minimizing $E_\bbA$ versus $E_\rA$ is presented in Section \ref{sec:6.1}.

\begin{figure}[tbp]
\centering
\includegraphics{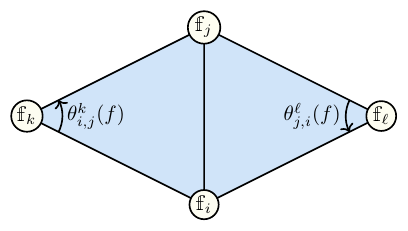}
\caption{An illustration of the cotangent weight defined on the surface $f(\M)$.}
\label{fig:cot}
\end{figure}

In practice, a spherical simplicial mapping $f$ is represented by an $n\times 3$ matrix
\begin{equation*}
\bf \equiv 
\begin{bmatrix}
    f^1_1   &   f^2_1  &  f^3_1\\
    \vdots     &   \vdots    &  \vdots\\
    f^1_n   &   f^2_n  &  f^3_n
\end{bmatrix} \equiv \big[ \bf^1, \bf^2, \bf^3 \big] \equiv
\begin{bmatrix}
    \bf_1^\top \\
    \vdots\\
    \bf_n^\top
\end{bmatrix},
\label{eq:f matrix}
\end{equation*}
where $\bf_i \in \mathbb S^2$ for $i = 1, \ldots, n$. 
From \cite[Lemma 3.1]{Yueh23}, the stretch energy $E_\rS$ can be expressed as
\[
E_\rS(\bf) = \frac{1}{2} \sum_{s = 1}^3  {\bf^s}^\top L_S(\bf) \,\bf^s,
\]
where $L_S(\bf)$ is the weighted Laplacian matrix given by
\begin{equation*} \label{eq:Ls}
[L_S(\bf)]_{i, j}=
\begin{cases}
    -\frac{1}{2} \big(\frac{\cot(\theta_{i, j}^k(\bf))}{\sigma_{\bf}([v_i, v_j, v_k])} + \frac{\cot (\theta_{j, i}^\ell (\bf))}{\sigma_{\bf}([v_j, v_i, v_\ell])} \big)   &\text{if}~ [v_i, v_j] \in \mathcal{E(M)},\\
    -\sum_{\ell \neq i} [L_S(\bf)]_{i, \ell}     &\text{if}~j=i,\\
    0       &\text{otherwise,}
\end{cases}
\end{equation*}
in which $\theta_{i,j}^k (\bf)$ is the angle opposite to the edge $[\bf_i, \bf_j]$ at the point $\bf_k$, as illustrated in Figure~\ref{fig:cot}, and the stretch factor $\sigma_{\bf}([v_i, v_j, v_k])$ is expressed as
$$
\sigma_{\bf}([v_i, v_j, v_k]) = \frac{|[v_i, v_j, v_k]|}{|[\bf_i, \bf_j, \bf_k]|}.
$$
As a result, the minimization of spherical authalic energy \eqref{eq:spherical Ea} is formulated as 
\begin{equation}
\argmin_{ \substack{\bf_i\in \mathbb{S}^2\\ i = 1,\ldots, n} } E_\bbA(\bf) \equiv\frac{|\M|}{6\, \bV(\bf)} \sum_{s = 1}^3 \bf^s L_S(\bf) \, \bf^s - 3 \, \bV(\bf).
\label{eq:obj_fun}
\end{equation}

Note that the constraint of $\bf_i\in \mathbb{S}^2$ for $i = 1,\ldots, n$ can alternatively be expressed using spherical coordinates, as discussed in the next section.

\subsection{Reformulation in spherical coordinates}

To solve the constraint optimization in \eqref{eq:obj_fun}, we reduce it to an unconstrained problem by globally parameterizing each $\bf_i$ using spherical coordinates $(\theta_i, \phi_i)$. This reformulation allows for unconstrained optimization, but requires computing gradients with respect to $\theta_i$ and $\phi_i$ instead of $\bf_i$.

Specifically, $\bf$ can be represented as 
\begin{equation}
\bf^1 = \sin \btheta \odot \cos \bphi, ~~ \bf^2 = \sin \btheta \odot \sin \bphi, ~~\text{and}~ \bf^3 = \cos \btheta,
\label{eq:sphere_coor}
\end{equation}
where $\odot$ denotes the { {\it Hadamard product} (or {\it Schur product}) \cite{Davi62} of vectors defined as
$$
(a_1, \ldots, a_n)^\top \odot (b_1, \ldots, b_n)^\top = (a_1 b_1, \ldots, a_n b_n)^\top,
$$}
$\btheta = (\theta_1, \ldots, \theta_{n})^{\top}$, and $\bphi = (\phi_1, \ldots, \phi_{n})^{\top}$, with $\theta_i \in [0, \pi]$, and $\bphi \in (-\pi, \pi]$, for $i = 1, \ldots, n$.
The inverse transformation is given by
\begin{equation}
\btheta = \arccos(\bf^3), ~~ \bphi = \mathrm{atan2}(\bf^2, \bf^1),
\label{eq:sphere_coor_inv}
\end{equation}
where $\mathrm{atan2}$ denotes the four-quadrant inverse tangent.
Therefore, $\bf$ can be represented in a vector as  
\begin{equation}
\f = 
\begin{bmatrix}
    \btheta \\
    \bphi
\end{bmatrix} \in \mathbb R^{2n}.
\label{eq:sphere_coor f}
\end{equation}

\begin{figure}[tbp]
\centering
\includegraphics{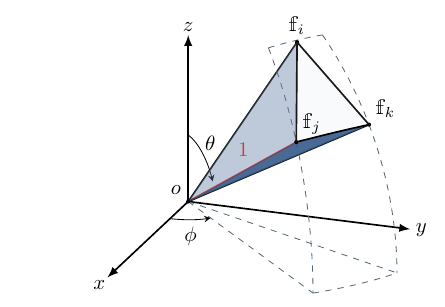}
\caption{An illustration for the tetrahedron $[o, \bf_i, \bf_j, \bf_k]$ formed by the face $[\bf_i, \bf_j, \bf_k]$ and the origin $o$ of $\R^3$.}
\label{fig:1}
\end{figure}

The volume measure $\bV$, as defined in \eqref{eq:volume measurement}, can be formulated as 
\begin{equation}
\bV(\bf) 
= \sum_{[v_i, v_j, v_k]\in\F(\M)} \frac{\bf_i^\top (\bf_j \times \bf_k)}{6},
\label{eq:Vol_formulate}
\end{equation}
as illustrated in Figure~\ref{fig:1}.
To compute the gradient of $|[o, \bf_i, \bf_j, \bf_k]|$, we let $\tau = \{i, j, k\}$ and denote the $x$, $y$, and $z$-coordinate of $\bf_i$, $\bf_j$, and $\bf_k$ by
\begin{subequations} \label{eq:4:16}
\begin{align}
\bf_\tau^1 &= 
\begin{bmatrix}
\sin \theta_i \cos \phi_i \\
\sin \theta_j \cos \phi_j \\
\sin \theta_k \cos \phi_k \\
\end{bmatrix} 
= \begin{bmatrix}
\sin \theta_i \\
\sin \theta_j \\
\sin \theta_k \\
\end{bmatrix}
\odot
\begin{bmatrix}
\cos \phi_i \\
\cos \phi_j \\
\cos \phi_k \\
\end{bmatrix}
= \sin\btheta_\tau \odot \cos\bphi_\tau,\\
\bf_\tau^2 &= 
\begin{bmatrix}
\sin \theta_i \sin \phi_i \\
\sin \theta_j \sin \phi_j \\
\sin \theta_k \sin \phi_k \\
\end{bmatrix} 
= \begin{bmatrix}
\sin \theta_i \\
\sin \theta_j \\
\sin \theta_k \\
\end{bmatrix}
\odot
\begin{bmatrix}
\sin \phi_i \\
\sin \phi_j \\
\sin \phi_k \\
\end{bmatrix}
=  \sin\btheta_\tau \odot \sin\bphi_\tau,\\
\bf_\tau^3 &= 
\begin{bmatrix}
\cos \theta_i \\ 
\cos \theta_j \\ 
\cos \theta_k 
\end{bmatrix} = \cos \btheta_\tau.
\end{align}
\end{subequations}
From \eqref{eq:Vol_formulate}, the gradient of $|[o, \bf_i, \bf_j, \bf_k]|$ with respect to $\bf_\tau^1$, $\bf_\tau^2$, and $\bf_\tau^3$ are given by
\begin{subequations} \label{eq:4:17}
\begin{align}
   \nabla_{\bf_\tau^1} |[o, \bf_i, \bf_j, \bf_k] | &= \frac{1}{6}(\bf_\tau^2 \times \bf_\tau^3),\\
   \nabla_{\bf_\tau^2} |[o, \bf_i, \bf_j, \bf_k] | &= \frac{1}{6}(\bf_\tau^3 \times \bf_\tau^1),\\
   \nabla_{\bf_\tau^3} |[o, \bf_i, \bf_j, \bf_k] | &= \frac{1}{6}(\bf_\tau^1 \times \bf_\tau^2).
\end{align}
\end{subequations}
Then, by applying the chain rule with \eqref{eq:4:16} and \eqref{eq:4:17}, the gradients of $|[o, \bf_i, \bf_j, \bf_k] |$ with respect to $\btheta_\tau$ and $\bphi_\tau$ are formulated as
\begin{subequations} \label{eq:grad_vol}
\begin{align}
\nabla_{\btheta_\tau} |[o, \bf_i, \bf_j, \bf_k] | 
&= \frac{1}{6} \big( \cos\btheta_\tau \odot \cos\bphi_\tau \odot (\bf^2_\tau \times \bf^3_\tau) \nonumber \\
&~~~~\, + \cos\btheta_\tau \odot \sin\bphi_\tau \odot (\bf^3_\tau \times \bf^1_\tau)
- \sin\btheta_\tau \odot (\bf^1_\tau \times \bf^2_\tau) \big),
\end{align}
and
\begin{align}
\nabla_{\bphi_\tau} |[o, \bf_i, \bf_j, \bf_k] | 
&= \frac{1}{6} \big( \sin\btheta_\tau \odot \cos\bphi_\tau \odot (\bf^3_\tau \times \bf^1_\tau) 
- \sin\btheta_\tau \odot \sin\bphi_\tau \odot (\bf^2_\tau \times \bf^3_\tau) \big).
\end{align}
\end{subequations}
Consequently, by assembling the contributions from all relevant tetrahedra using the formulation in \eqref{eq:grad_vol}, we can compute the gradient of the volume $\bV(\bf)$. Denote $\bV(\f)$ the volume expressed in the spherical coordinates $\f$ from \eqref{eq:sphere_coor f}. We then write $\nabla_{\btheta} \bV(\f)$ and $\nabla_{\bphi} \bV(\f)$ for the gradients of $\bV(\f)$ with respect to $\btheta$ and $\bphi$, respectively.

From \cite[Theorem 3.5]{Yueh23}, the gradient of $E_\rS(\bf)$ with respect to $\bf$ is given by
\begin{equation} 
\nabla_{\bf^s} E_\rS(\bf) = 2 \, L_S(\bf) \,\bf^s, ~~~ \text{for $s=1,2,3$}.
\label{eq:GradEs}
\end{equation}

By applying the chain rule along with \eqref{eq:4:16} and \eqref{eq:GradEs}, the gradients of $E_\bbA$ with respect to $\f$ are formulated as
\begin{subequations} \label{eq:Grad}
\begin{align}
\nabla_{\btheta} E_\bbA (\f) 
&= \Big( \frac{|\M|}{3\bV(\f)} \Big) \nabla_{\btheta} E_\rS(\f) + \nabla_{\btheta} \Big( \frac{|\M|}{3\bV(\f)} \Big) E_\rS(\f) - \nabla_{\btheta} \big(3\bV(\f) \big) \nonumber\\
&= \Big( \frac{2\, |\M|}{3\bV(\f)} \Big) \Big( \cos\btheta \odot \cos\bphi \odot L_S(\f) \, \bf^1  + \cos\btheta \odot \sin\bphi \odot L_S(\f) \, \bf^2 \nonumber\\
&~~~ -\sin \btheta \odot L_S(\f) \, \bf^3 \Big) - \Big( 3+ \frac{|\M|\, E_\rS(\f)}{3\bV(\f)^2} \Big) \nabla_\btheta \bV(\f). 
\end{align}
and
\begin{align}
\nabla_{\bphi} E_\bbA (\f) 
&= \Big( \frac{|\M|}{3\bV(\f)} \Big) \nabla_{\bphi} E_\rS(\f) + \nabla_{\bphi} \Big( \frac{|\M|}{3\bV(\f)} \Big) E_\rS(\f) - \nabla_{\bphi} 3\bV(\f) \nonumber\\
&= \Big( \frac{2\, |\M|}{3\bV(\f)} \Big) \Big( \sin\btheta \odot \cos\bphi \odot L_S(\f) \, \bf^2  - \sin\btheta \odot \sin\bphi \odot L_S(\f) \, \bf^1 \Big)\nonumber\\
 &~~~ - \Big( 3+ \frac{|\M|\, E_\rS(\f)}{3\bV(\f)^2} \Big) \nabla_\btheta \bV(\f). 
\end{align}
\end{subequations}

Using the explicit gradient formula of $E_\bbA$ given in \eqref{eq:Grad}, we develop an associated algorithm to solve the optimization problem \eqref{eq:obj_fun}. The details of this method are presented in the following section.

\section{Preconditioned nonlinear CG method}
\label{sec:4}

The { \textit{conjugate gradient (CG) method}}, introduced by Hestenes and Stiefel \cite{HeSt52}, was originally developed to solve large-scale linear systems with a positive-definite coefficient matrix. Fletcher and Reeves \cite{FlRe64} later adapted it for nonlinear optimization. In this section, we present the algorithmic procedure of the nonlinear CG method with suitable preconditioning to solve the objective problem \eqref{eq:obj_fun}, and then discuss its convergence.

\subsection{Algorithmic procedure}
\label{sec:4.1}

{ 
We first compute the initial mapping using the fixed-point method incorporating stereographic projection \cite[Algorithm 4.3]{YuLL19} for a few iterations. While the method effectively reduces the energy at the first few iterations, its performance diminishes rapidly and may occasionally lead to an increase in energy. To address this issue, we limit the algorithm \cite[Algorithm 4.3]{YuLL19} to at most 15 iterations or terminate it early if the energy $E_\rA$ increases. A comparison of numerical results under different initial mappings is presented in Section~\ref{sec:6.2}.
}

To eliminate rotational degrees of freedom, we fix two vertices during the optimization. Specifically, for each vertex, we compute the ratio $|f(\tau)| / |\tau|$ over its 1-ring and then select the two vertices whose 1-ring ratios are closest to the global mean. For simplicity, the spherical coordinates of the remaining vertices are still denoted by $\f$ as in \eqref{eq:sphere_coor f}.

The CG method is a { \em line search method} \cite[Chapter 3]{NoWr06}, meaning that at each iteration $k \geq 0$, the vector $\f^{(k)}$ is updated via a step length $\alpha_k \in \mathbb{R}$ along a search direction $\p^{(k)} \in \mathbb{R}^{2(n-2)}$:
\begin{equation}  \label{eq:update}
    \f^{(k+1)} = \f^{(k)} + \alpha_k \p^{(k)}. 
\end{equation}

The initial direction of our preconditioned nonlinear CG method is the preconditioned negative gradient
\begin{equation} \label{eq:preconditioner}
    \p^{(0)} = -M^{-1} \g^{(0)}, ~~~ M = I_2 \otimes L_S(\f^{(0)}),
\end{equation}
where $\g^{(0)}$ is the gradient $ \nabla_{\f} E_\bbA(\f^{(0)})$ defined in \eqref{eq:Grad}.
{ 

Note that the preconditioner $M$, constructed using the submatrix of the Laplacian $L_S(\f^{(0)})$ by fixing two vertices, is typically symmetric positive definite. This property ensures that the initial search direction $\p^{(0)}$ is a descent direction. 
}
For the $k$th iteration, $k \geq 1$, we add a correction term based on the previous direction $\p^{(k-1)}$ scaled by the coefficient $\beta_k$:
\begin{equation} \label{eq:direction}
    \p^{(k)} = -M^{-1} \g^{(k)} + \beta_k \p^{(k-1)}, ~~~ \beta_k = \frac{{\g^{(k)}}^\top M^{-1} \g^{(k)} }{ {\g^{(k-1)}}^\top M^{-1} \g^{(k-1)} }.
\end{equation}

{ 
To efficiently solve the linear systems in practice, we perform a {\it reordered Cholesky decomposition} \cite{ChDH08} of $L_S(\f^{(0)})$, given by
\begin{equation} \label{eq:LS3}
U^\top U = P^\top L_S(\f^{(0)}) P,
\end{equation}
where $U$ is an upper triangular matrix and $P$ is a permutation matrix obtained via the {\it approximate minimum degree} ordering \cite{AmDD04}.
This factorization allows us to efficiently solve systems of the form $L_S(\f^{(0)}) \, \x = \r$ by solving the triangular systems:
\begin{subequations} \label{eq:LS}
\begin{equation} \label{eq:LS1}
U^\top \y = P^\top \r \,\text{ and }\, U\z = \y,
\end{equation}
followed by recovering the solution via
\begin{equation} \label{eq:LS2}
\x = P\z. 
\end{equation}
\end{subequations}
} 

For the step length $\alpha_k$, we consider the 1-dimensional function
\[
\varphi (\alpha) =  E_\bbA(\f^{(k)} + \alpha \p^{(k)}).
\]
Ideally, the optimal step length is $\alpha = \argmin \varphi(\alpha)$. However, due to the complexity of $E_\rS$, the minimizer does not have an explicit form.  Hence, we approximate $\varphi(\alpha)$ by interpolating a quadratic polynomial and use its minimizer as an estimate for $ argmin \varphi(\alpha)$ \cite[Section 3.5]{NoWr06}.

Specifically, using the initial guess $\alpha_{k-1}$, we construct a quadratic function $\zeta(\alpha)$ by interpolating the three available information:
\begin{subequations} \label{eq:Step}
\begin{equation} \label{eq:Phi_cond}
\begin{cases}
\varphi(0) =  E_\bbA ({\f}^{(k)}), \\
\varphi'(0) = \tfrac{\d}{\d\alpha} E_\bbA ({\f}^{(k)} + \alpha \p^{(k)})\big|_{\alpha = 0} = {{\p}^{(k)}}^\top \g^{(k)}, \\
\varphi(\alpha_{k-1}) = E_\bbA ({\f}^{(k)} + \alpha^{(k-1)} \p^{(k)}).
\end{cases}
\end{equation}
We then obtain
\begin{equation}
\zeta(\alpha) = a \alpha^2 + b \alpha + c \equiv 
\Big( \dfrac{\varphi(\alpha_{k-1}) - \varphi(0) - \alpha_{k-1} \varphi'(0) }{\alpha_{k-1}^2} \Big)\, \alpha^2 + \varphi'(0) \, \alpha + \varphi(0).
\end{equation}
Consequently, $\alpha_k$ is estimated the minimizer of $\zeta(\alpha)$, obtained by setting $\zeta'(\alpha) = 0$:
\begin{equation}
\alpha_k = - \frac{b}{2a}.
\end{equation}
\end{subequations}
If $\alpha_k$ does not sufficiently decrease the energy, we use the minimizer of $\zeta(\alpha)$ as a new initial guess and repeat the interpolation. The pseudo-code for the proposed method is summarized in Algorithm \ref{alg:PCG}.

\begin{algorithm}[ht]
\caption{Preconditioned nonlinear CG method for spherical authalic mapping}
\label{alg:PCG}
\begin{algorithmic}[1]
\Require A simplicial surface $\M$.
\Ensure A spherical authalic mapping $\bf^*$.
\State Perform \cite[Algorithm 4.3]{YuLL19} to 15 iterations.
\State Compute spherical coordinates $[\btheta, \bphi]$ by \eqref{eq:sphere_coor_inv}.
\State Compute gradient $\g$ by \eqref{eq:Grad} with respect to $\f$.
\State Compute the preconditioner $M$ by \eqref{eq:preconditioner}.
\State Perform Cholesky decompositions of $M$ by \eqref{eq:LS3}.
\State Solve $M\h = \g$ by \eqref{eq:LS}.
\State Let $\p = -\h$.
\State Let $\alpha = 0.01$.
\While {not converge}
    \State Update $\alpha$ as \eqref{eq:Step}.
    \State Update $\f \leftarrow \f + \alpha \p$.
    \State Let $\gamma = \h^\top \g$.
    \State Update the gradient $\g$ \eqref{eq:Grad}.
    \State Solve $M\h = \g$ by \eqref{eq:LS}.
    \State Update $\beta = (\h^{\top} \g) / \gamma$.
    \State Update $\p \gets -\h + \beta \p$.
\EndWhile
\State Obtain $\bf^*$ from $\f$ by \eqref{eq:sphere_coor}.
\end{algorithmic}
\end{algorithm}

\subsection{Global convergence}
The global convergence of Algorithm \ref{alg:PCG} is theoretically guaranteed under the assumption that each step length $\alpha_k$ is chosen to satisfy the {{\it strong Wolfe conditions} \cite[Chapter 3.1]{NoWr06}}:
\begin{subequations} \label{eq:Wolfe}
\begin{align}
{ E_\bbA(\f^{(k+1)}) - E_\bbA (\f^{(k)}) }& \leq  c_1 \, \alpha_k {\g^{(k)}}^{\top} \p^{(k)}, \label{eq:Wolfe1} \\
|{\g^{(k+1)}}^{\top}  \p^{(k)}| & \leq c_2 \, | {\g^{(k)}}^{\top} \p^{(k)}|, \label{eq:Wolfe2}
\end{align}
\end{subequations}
where $0 < c_1 < c_2 < \frac{1}{2}$, and $\g^{(k+1)}$ is the gradient at $\f^{(k+1)}$. 
Under this assumption, we establish the following result:

\begin{theorem} \label{thm:convergence}

Suppose Algorithm \ref{alg:PCG} (preconditioned nonlinear CG method for $E_\bbA$) is implemented with each step length satisfying the strong Wolfe conditions \eqref{eq:Wolfe} with $0<c_1<c_2<\frac{1}{2}$. Assume that the iterates $\{\f^{(k)}\}$ remain in a compact subset of the nondegenerate region. Then,
\[
\liminf_{k \rightarrow \infty} \| \g^{(k)} \|_{M^{-1}} = 0,
\]
{ 
provided $M$ is a symmetric positive definite matrix.
}
\end{theorem}
\begin{proof}
{  The nondegeneracy assumption keeps all image-face areas and the volume denominator $\bV$ uniformly away from zero. Hence, $E_\bbA$ is smooth on a neighborhood of the compact set containing the iterates, and both $E_\bbA$ and its gradient are Lipschitz continuous there. Moreover, $E_\bbA$ is bounded below by zero, and the preconditioner $M$ is a symmetric positive definite matrix (constructed from a submatrix of the Laplacian matrix). Therefore, the convergence result follows directly from \cite[Appendix A]{LiYu24}.}

\end{proof}

Theorem \ref{thm:convergence} shows that there is a subsequence $\{\g^{(k_\ell)}\}_{\ell=1}^\infty$ such that
$$
\lim_{\ell\to\infty} \|\g^{(k_\ell)}\|_{M^{-1}} = 0.
$$
Consequently, for any tolerance $\varepsilon>0$, one can find an index $\ell$ with $\| \g^{(k_\ell)} \|_{M^{-1}}<\varepsilon$. In particular, if the stopping rule is
$
\|\g^{(k)}\|_{M^{-1}} < \varepsilon,
$
then some $k_\ell$ will eventually meet it. Hence, Algorithm \ref{alg:PCG} must terminate for any prescribed tolerance.

\section{Riemannian bijective correction}
\label{sec:5}
In Section \ref{sec:3.1}, we { modified} the authalic energy $E_\rA$ in \eqref{eq:Ea} to the spherical authalic energy $E_\bbA$ in \eqref{eq:spherical Ea}. This crucial modification significantly improves bijectivity, as confirmed by numerical comparisons in Section \ref{sec:6.1}. 
To guarantee bijectivity, we introduce a { \em Riemannian bijective correction} as a post-processing step, a local correction method by solving a linear system on the tangent plane $T_{\c(\tau)}\mathbb{S}$ for each folded face $\tau$.

\begin{figure}[tbp]
\centering
\includegraphics{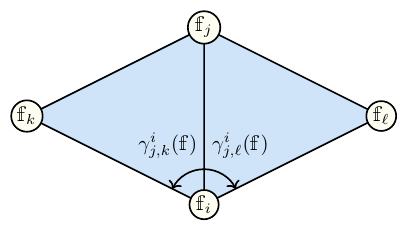}
\caption{An illustration of the mean-value weight \cite{Floa03b} defined on $f(\M)$.}
\label{fig:mean value}
\end{figure}

We firstly construct the Laplacian matrix using the { {\it mean-value weights}} \cite{Floa03b}, defined by
\begin{equation}
    [L_\mathrm{M}(\bf)]_{i,j} = 
    \begin{cases}
        - \sum_{[v_i, v_j, v_k] \in \F(\M)} \frac{\tan (\gamma_{j,k}^i(\bf) /2)}{\|\bf_i - \bf_j\|_2}  & \mbox{if}~[v_i, v_j] \in \E(\M),\\
        -\sum_{\ell \neq i} [L_\mathrm{M}(\bf)]_{i, \ell}  & \mbox{if}~j = i,\\
        0   &\mbox{otherwise},
    \end{cases}
\label{eq:L_M}
\end{equation}
where $\gamma_{j, k}^i(\bf)$ is the angle opposite to the edge $[\bf_j, \bf_k]$ at the vertex $\bf_i$ on the surface $f(\M)$ (see Figure~\ref{fig:mean value}).

Next, to correct a single folded triangle $f(\tau) = [\bf_i, \bf_j, \bf_k]$, let $\mathtt{F} = \{i, j, k\}$ denote the indices of its vertices and $\bF^c = \{1, \cdots, n\} \setminus \bF$ denote the complement. We denote the face center by $\c(\tau)$, defined in \eqref{eq:FaceCenter}. The vertices are then projected onto the tangent plane $T_{\c(\tau)}\mathbb{S}$ via

\begin{equation}
\widetilde \bf = \bh - \big( \bh \, \c(\tau) \big) \, \c(\tau)^\top + \1 \c(\tau)^\top, ~~~ \bh = \bf - \1 \c(\tau)^\top,
\end{equation}
where $\1 =  (1, \cdots, 1)^\top \in \mathbb R^n$. 

Let $[L_\mathrm{M}(\bf)]_{\bF, \bF^c}$ be the submatrix of $L_\mathrm{M}(\bf)$ with entries $[L_\mathrm{M}(\bf)]_{i,j}$ for $i \in \mathtt{F}$ and $j \in \mathtt{F}^c$, and let $\bf^s_\bF$ denote the subvector of $\bf$ composed of the $s$-th coordinates of vertices in $\mathtt{F}$. 
We unfold the folded triangle by solving the following 3-by-3 linear system for each coordinate $s = 1, 2, 3$:

\begin{equation}
    [L_\mathrm{M}(\bf)]_{\bF, \bF} \widetilde \bf^s_{\bF} = - [L_\mathrm{M}(\bf)]_{\bF, \bF^c} \widetilde \bf^s_{\bF^c}.
\label{eq:MVT_correct}
\end{equation}
The updated points $\widetilde \bf_{\bF}$ are then projected back onto the sphere:
\begin{equation*}
    \bf_i = \frac{\widetilde{\bf}_i}{\|\widetilde{\bf}_i\|_2}, \quad i \in \bF.
\end{equation*}

The correction step in \eqref{eq:MVT_correct} only depends on the 1-ring neighborhood of the folded face.

{
\begin{remark}
The term Riemannian bijective correction refers to the extension of the correction method on $\mathbb{S}^2$, which is a Riemannian $2$-manifold, rather than in Euclidean space.
\end{remark}
}

This method ensures the unfolding of locally folded faces, similarly to the planar case for open surfaces discussed in \cite{Yueh23}. We formalize this in the following theorem:
\begin{theorem} \label{thm:bij_cor}
Suppose the projection of the 1-ring of a face onto the tangent plane is convex. Then, the updated spherical mapping is locally bijective.

\end{theorem}
\begin{proof}
Let $\bF$ be a folded face. It suffices to show that the 1-ring of the updated mapping $\bf_{\bF}$ is bijective. 
We denote $\bN_i=\{ j \mid [v_i,v_j]\in\E(\M) \}$ as the vertex indices of the 1-ring neighborhood of the vertex $v_i\in\V(\M)$. 
From \eqref{eq:MVT_correct}, we have $[L_\mathrm{M}(\bf)\, \widetilde \bf]_\bF = \0$, and then for each $i\in\bF$,
$$
\widetilde \bf_i = \sum_{\ell\in\bN_i} \frac{-[L_\mathrm{M}(\bf)]_{i,\ell}}{[L_\mathrm{M}(\bf)]_{i,i}} \, \widetilde \bf_\ell.
$$
The inequality $0 < \gamma_{j,k}^i (\bf) < \pi$ implies that $\tan (\gamma_{j,k}^i(\bf) /2) > 0$. Therefore, the mean value weight in \eqref{eq:L_M} is always positive so that $-[L_\mathrm{M}(\bf)]_{i,\ell}>0$ for $\ell\neq i$. 
This implies that $\frac{-[L_\mathrm{M}(\bf)]_{i,\ell}}{[L_\mathrm{M}(\bf)]_{i,i}}>0$ and $\sum_{\ell\in\bN_i}\frac{-[L_\mathrm{M}(\bf)]_{i,\ell}}{[L_\mathrm{M}(\bf)]_{i,i}}=1$.
Given the assumption that the projection of the 1-ring on the folded face $\bF$ is convex, it follows from \cite[Theorem 4.1]{Floa03} that the 1-ring of $\widetilde{\bf}_{\bF}$ is bijective. Note that the bijectivity of the projection on the tangent plane ensures the bijectivity of the spherical mapping, and hence the 1-ring of $\bf_{\bF}$ is bijective.
\end{proof}

In practice, the assumption may not always hold, therefore, this procedure is repeated iteratively for all folded triangles until the mapping becomes bijective. The detailed procedure is summarized in Algorithm \ref{alg:bij_cor}. Numerical experiments, demonstrated in Section \ref{sec:6.1}, show that the proposed iterative bijective correction successfully clears out all folded triangles for our non-bijective mappings.

\begin{algorithm}[t]
\caption{Riemannian bijective correction for spherical mappings}
\label{alg:bij_cor}
\begin{algorithmic}[1]
\Require A simplicial surface $\M$ and a folded spherical map $\bf$.
\Ensure A bijective spherical map $\bf^*$.
\State Compute the folded face's indices.
\While {$\#$ Foldings $>0$}
    \State Construct the mean value Laplacian $L_\mathrm{M}$ by \eqref{eq:L_M}. 
    \For {$k = 1, 2, \cdots, \#$ Foldings}
        \State Select one folded index $\mathtt{F}$.
        \State Set $\mathtt{F}^c = \{1, \cdots, n\} \setminus \mathtt{F}$.
        \State Compute the folding face center $\c$ as \eqref{eq:FaceCenter}.
        \State Set $\bh = \bf - \c^\top$.
        \State Update $\bh \leftarrow \bh - \big( \bh \,\c \big)\, \c^\top + \c^\top$.

        \State Solving the linear system $$ [L_\mathrm{M}]_{\bF, \bF} \bh^s_{\bF} = - [L_\mathrm{M}]_{\bF, \bF^c} \bh^s_{\bF^c}, ~~~\mbox{for}~s = 1,2,3.$$
        \State Update $\bf_\bF = \bh_{\bF} / \|\bh_{\bF}\|_2$.
    \EndFor
    \State Compute the folded face's indices.
\EndWhile
\State Obtain $\bf^*$ from $\bf$.
\end{algorithmic}
\end{algorithm}

\begin{figure}[p]
\centering
\resizebox{\textwidth}{!}{
\begin{tabular}{ccccccc}

\specialrule{.2em}{.1em}{.1em}
\multicolumn{3}{c}{Right Hand} && \multicolumn{3}{c}{David Head} \\
\multicolumn{3}{c}{$\F(\M) = 8,808$ $\V(\M) = 4,406$} && \multicolumn{3}{c}{$\F(\M) = 21,338$ $\V(\M) = 10,671$} \\
\includegraphics[height=3cm]{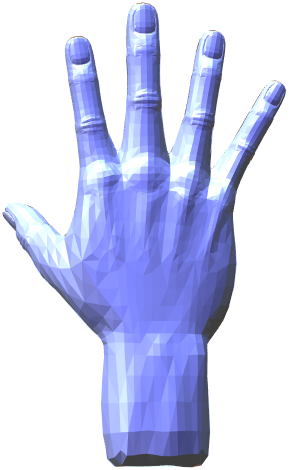} &
\includegraphics[height=3cm]{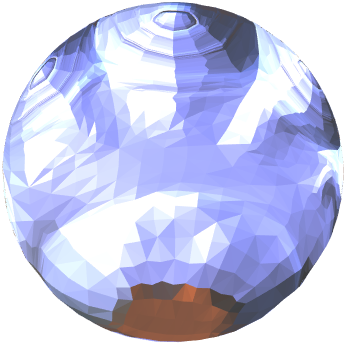} &
\includegraphics[height=3cm]{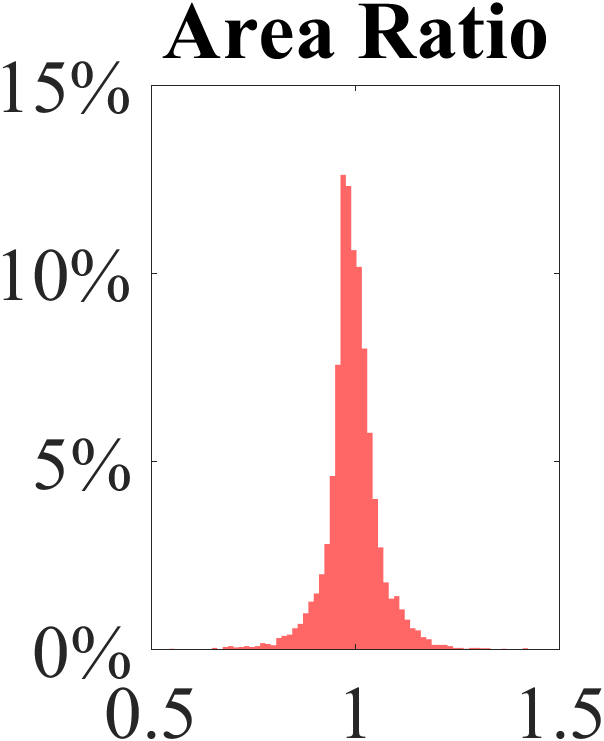} &&
\includegraphics[height=3cm]{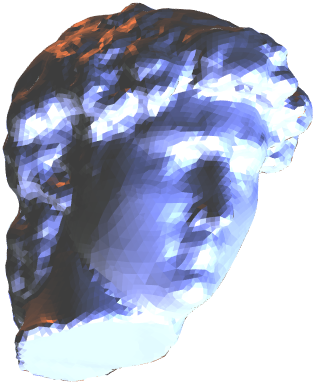} &
\includegraphics[height=3cm]{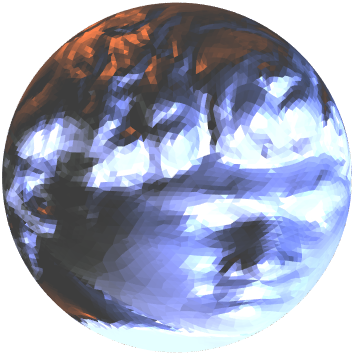} &
\includegraphics[height=3cm]{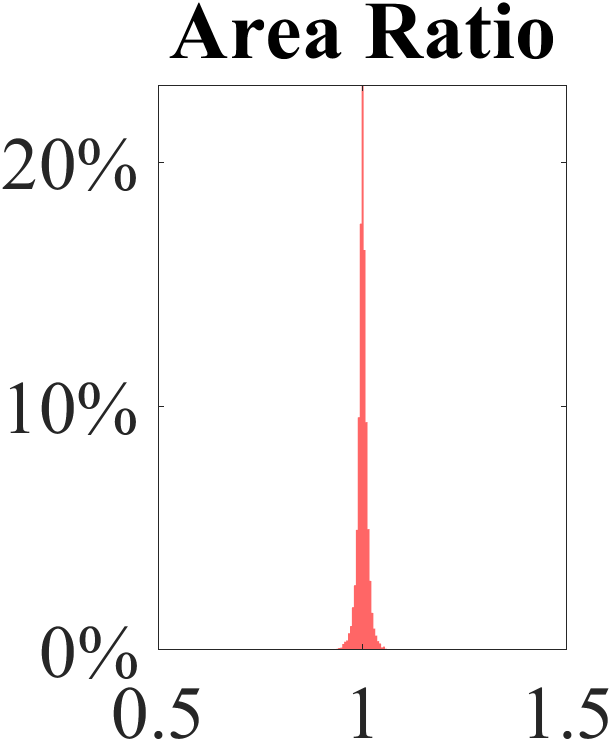} \\
\Cline{1pt}{1-3}\Cline{1pt}{5-7}
\multicolumn{3}{c}{Bull} && \multicolumn{3}{c}{Bulldog} \\
\multicolumn{3}{c}{$\F(\M) = 34,504$ $\V(\M) = 17,254$} && \multicolumn{3}{c}{$\F(\M) = 99,590$ $\V(\M) = 49,797$} \\
\includegraphics[height=3cm]{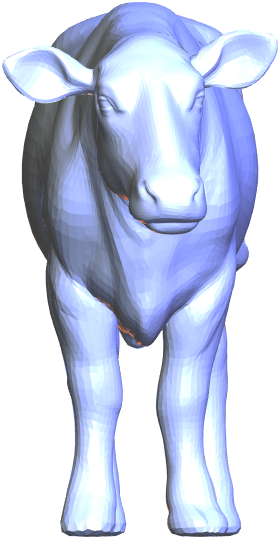} &
\includegraphics[height=3cm]{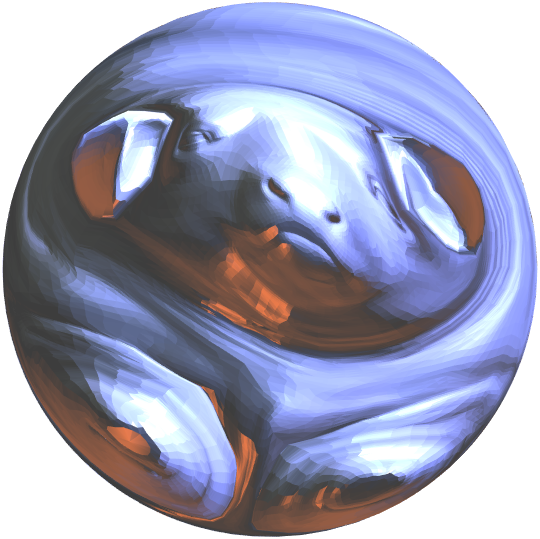} &
\includegraphics[height=3cm]{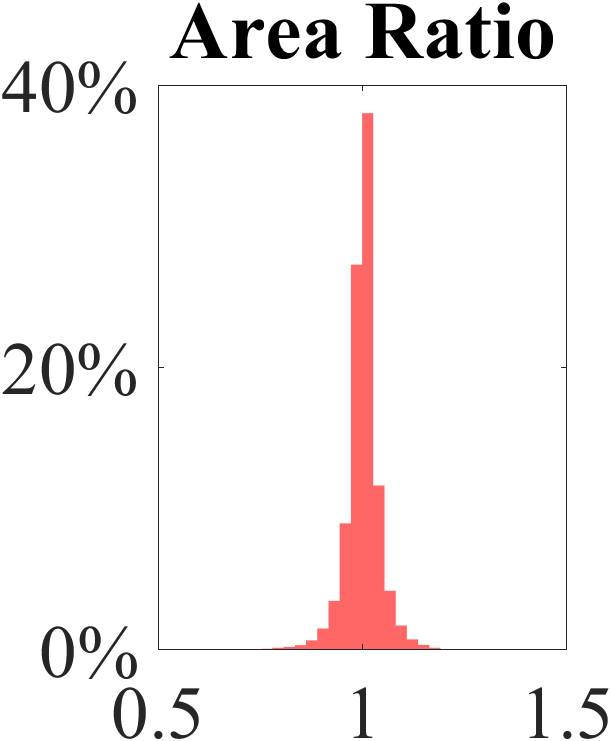} &&
\includegraphics[height=3cm]{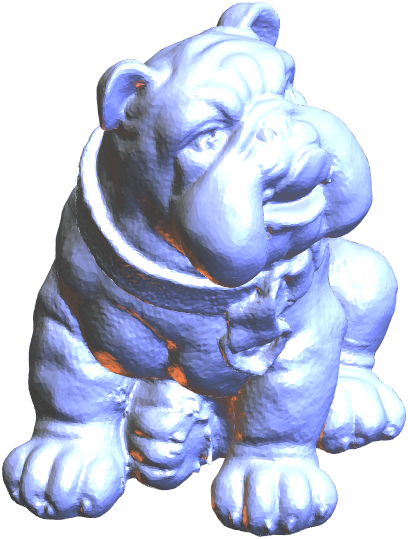} &
\includegraphics[height=3cm]{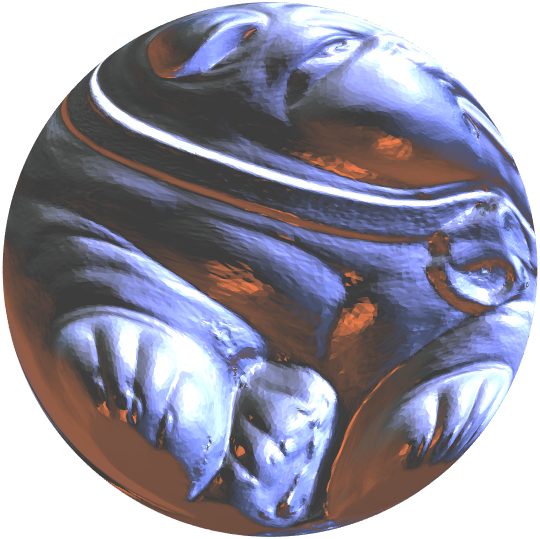} &
\includegraphics[height=3cm]{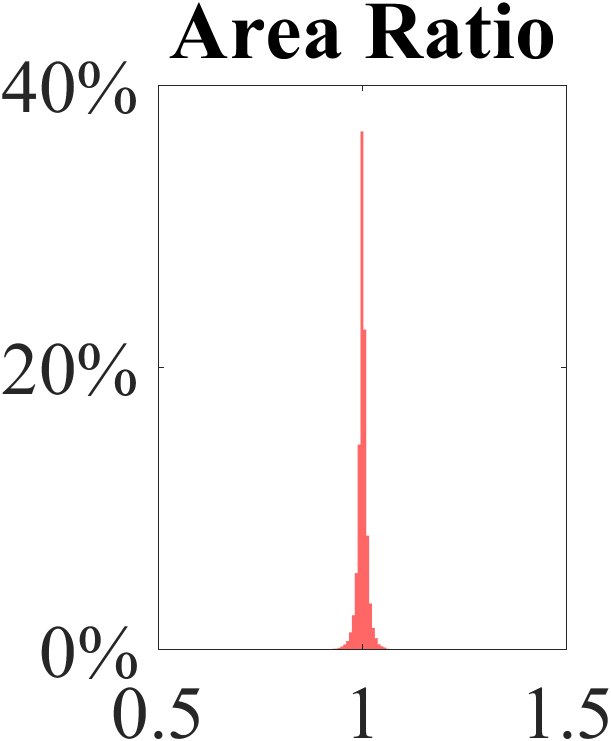} \\
\Cline{1pt}{1-3}\Cline{1pt}{5-7}
\multicolumn{3}{c}{Lion} && \multicolumn{3}{c}{Gargoyle} \\
\multicolumn{3}{c}{$\F(\M) = 100,000$ $\V(\M) = 50,002$} && \multicolumn{3}{c}{$\F(\M) = 100,000$ $\V(\M) = 50,002$} \\
\includegraphics[height=3cm]{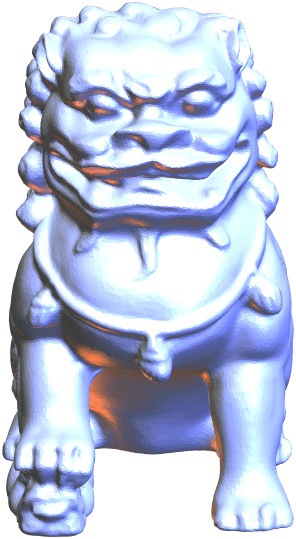} &
\includegraphics[height=3cm]{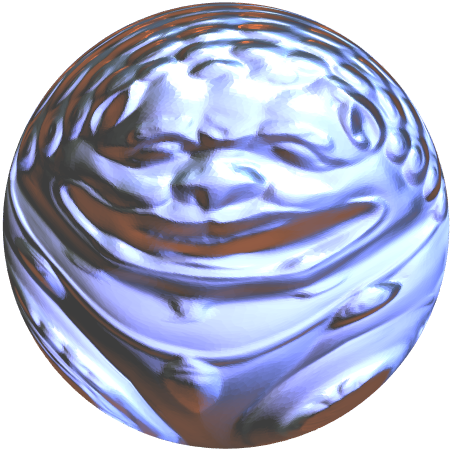} &
\includegraphics[height=3cm]{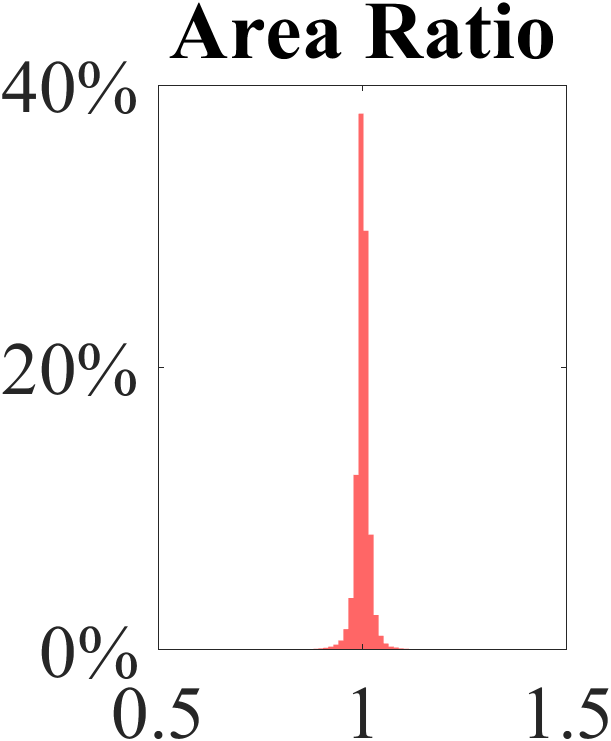} &&
\includegraphics[height=3cm]{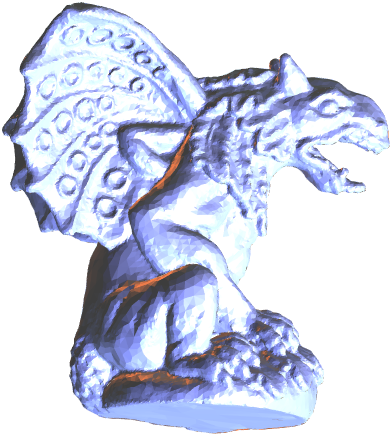} &
\includegraphics[height=3cm]{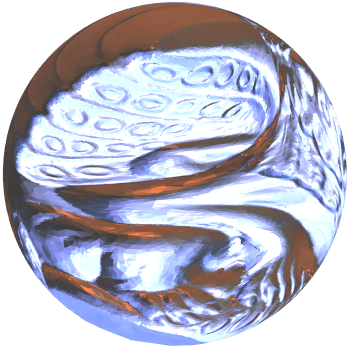} &
\includegraphics[height=3cm]{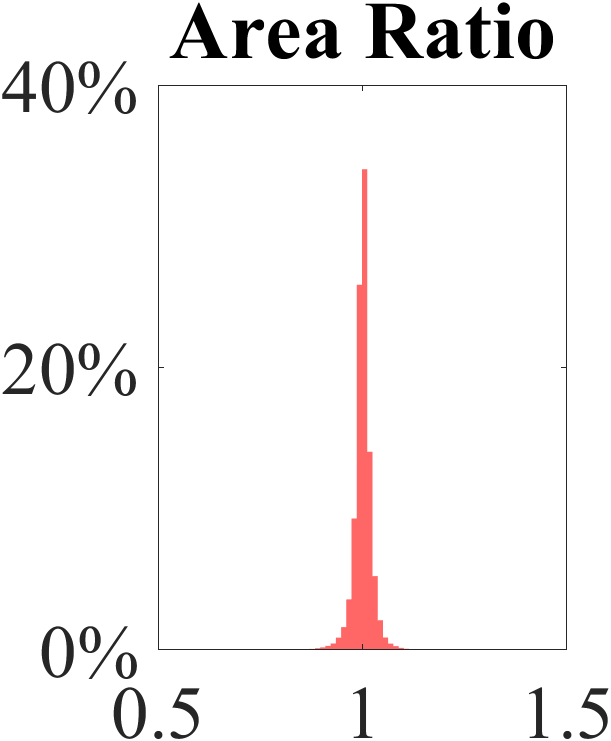} \\
\Cline{1pt}{1-3}\Cline{1pt}{5-7}
\multicolumn{3}{c}{Max Planck} && \multicolumn{3}{c}{Chess King} \\
\multicolumn{3}{c}{$\F(\M) = 102,212$ $\V(\M) = 51,108$} && \multicolumn{3}{c}{$\F(\M) = 263,712$ $\V(\M) = 131,858$} \\
\includegraphics[height=3cm]{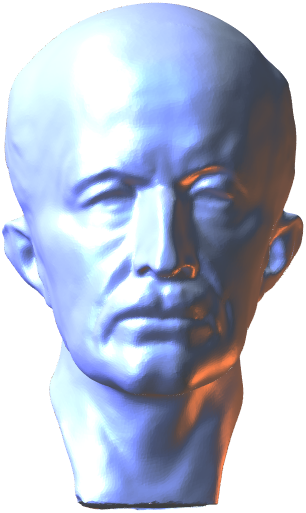} &
\includegraphics[height=3cm]{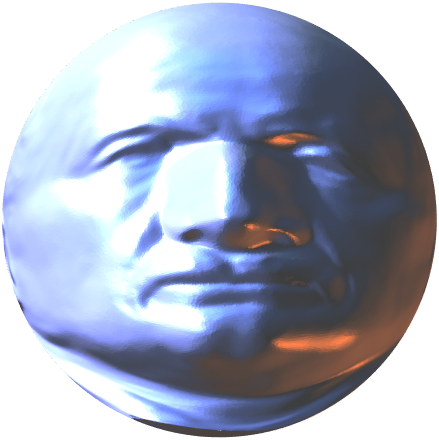} &
\includegraphics[height=3cm]{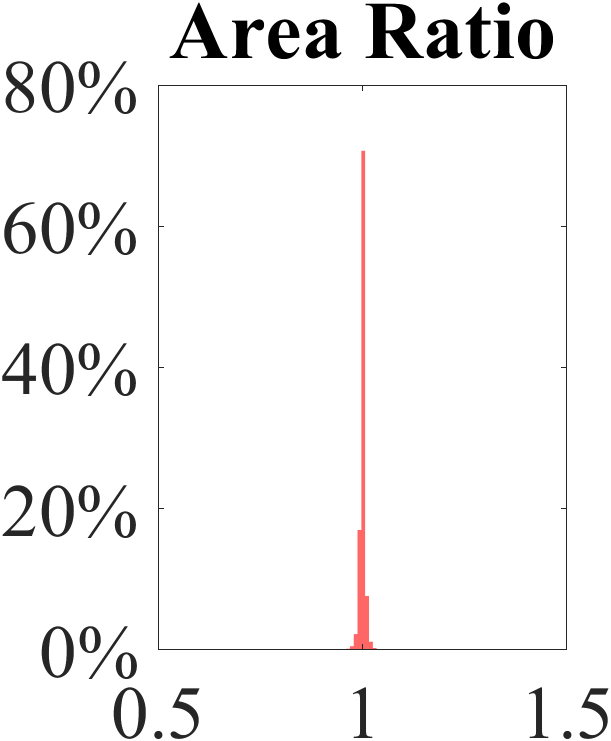} &&
\includegraphics[height=3cm]{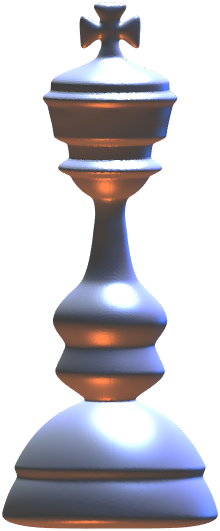} &
\includegraphics[height=3cm]{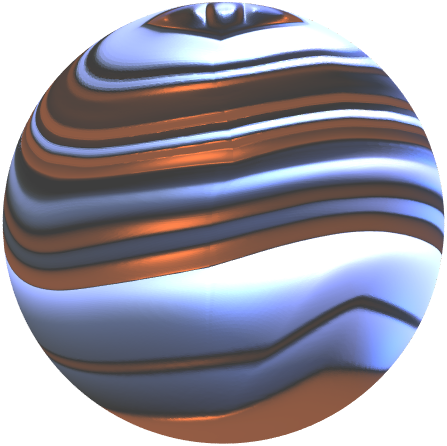} &
\includegraphics[height=3cm]{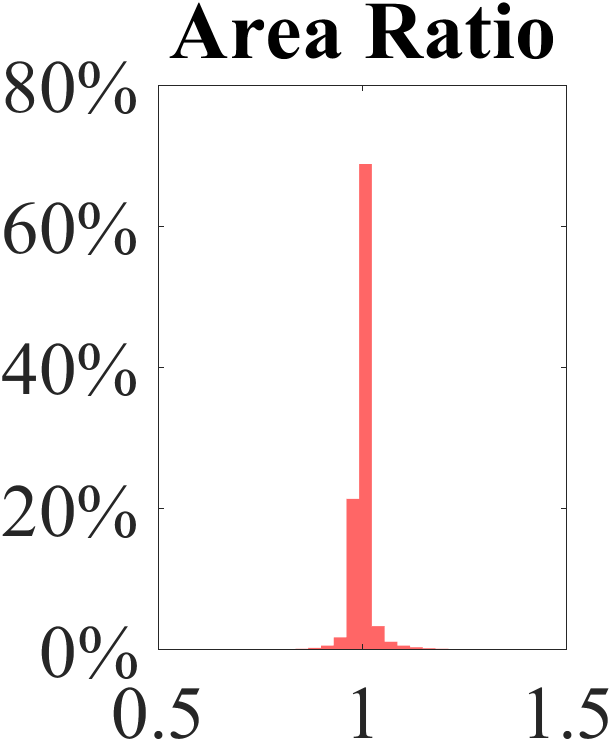} \\
\specialrule{.2em}{.1em}{.1em}
\end{tabular}
}
\caption{The benchmark mesh models along with their spherical area-preserving parameterization computed by our method, and histogram of their area ratio \eqref{eq:Area_dist}.}
\label{fig:MeshModel}
\end{figure}

\section{Numerical experiments}
\label{sec:6}
In this section, we present the numerical results of the preconditioned nonlinear CG method (Algorithm \ref{alg:PCG}) and the Riemannian bijective correction (Algorithm \ref{alg:bij_cor}) for spherical authalic (area-preserving) parameterizations.

Figure~\ref{fig:MeshModel} illustrates the benchmark mesh models along with their corresponding spherical authalic parameterizations. These models were obtained from well-known repositories, including the Stanford 3D Scanning Repository \cite{Stanford} and Sketchfab \cite{Sketchfab}. All experiments were conducted in MATLAB on a laptop equipped with an AMD Ryzen 9 5900HS processor and 32 GB of RAM.

In the following subsections, we quantify the local area distortion using the area ratio
\begin{equation}
    \mathcal{R}_\mathrm{area}(\tau) =  \bigg| \frac{|f(\tau)| / \A(f) }{|\tau| / |\M|} \bigg|.
\label{eq:Area_dist}
\end{equation}
For consistency across experiments, we normalize the area of the surface $\M$ to $4\pi$, ensuring that the mean of $\mathcal{R}_\mathrm{area}$ remains close to 1. Note that a perfectly area-preserving mapping satisfies
$$
\SD_{\tau \in \F(\M)} \mathcal{R}_\mathrm{area}(\tau) = 0.
$$

Figure~\ref{fig:MeshModel} shows the histogram of $\mathcal{R}_\mathrm{area}$ for the resulting parameterizations. We observe that these values are tightly concentrated around 1, indicating strong area preservation.
We also quantify the global area distortion via the authalic energy $E_\rA$ defined in \eqref{eq:Ea}. As discussed in Section \ref{sec:2.2}, $E_\rA(f) = 0$ if $f$ is perfectly area-preserving.

\subsection{Spherical authalic energy and bijective correction}
\label{sec:6.1}

{

As discussed in Section~\ref{sec:3.1}, we modify the original authalic energy $E_\rA$ to the spherical authalic energy $E_\bbA$ to improve bijectivity in the resulting mappings. We perform numerical experiments comparing the minimization of $E_\bbA$ and $E_\rA$ using Algorithm~\ref{alg:PCG}, with results summarized in Table~\ref{tab:Ea_comparison}. Notably, all mappings obtained by minimizing $E_\bbA$ are bijective, whereas minimizing $E_\rA$ sometimes produces hundreds of folded triangles. Moreover, minimizing $E_\bbA$ achieves slightly better area preservation.

We further evaluate the performance of the Riemannian bijective correction (Algorithm~\ref{alg:bij_cor}) as a post-processing step for non-bijective mappings obtained by minimizing $E_\rA$ on the Right Hand, Bull, Lion Statue, and Chess King mesh models. The numerical results are presented in Table~\ref{tab:bij_correct}. Remarkably, the correction successfully eliminates all folded triangles, including 165 in the Right Hand mesh model and 381 in the Chess King mesh model. However, this correction introduces a slight increase in area distortion, as reflected in both $E_\rA$ and the standard deviation of $\mathcal{R}_\mathrm{area}$. 

In summary, the results highlight the advantages of minimizing $E_\bbA$ and validate the Riemannian bijective correction as a reliable post-processing strategy.

}

\begin{table}[ht]
\centering
\caption{ {Numerical comparison of minimizing $E_\bbA$ \eqref{eq:spherical Ea} versus minimizing $E_\rA$ \eqref{eq:Ea} across benchmark mesh models.
$\mathcal{R}_\mathrm{area}$: area ratio \eqref{eq:Area_dist}; $\#$Foldings: number of folded triangles. }
}
\label{tab:Ea_comparison}
\resizebox{\textwidth}{!}{
\begin{tabular}{lrccrlrccr}
\specialrule{.2em}{.1em}{.1em}
\multirow{3}{*}{Model name} & \multicolumn{4}{c}{minimizing $E_\bbA$$^\dagger$}  & & \multicolumn{4}{c}{minimizing $E_\rA$$^\dagger$} \\
\cline{2-5} \cline{7-10}
& Time       & \multirow{2}{*}{$E_\rA$}  &$\mathcal{R}_\mathrm{area}$ &  $\#$Fold- & & Time       & \multirow{2}{*}{$E_\rA$}  &$\mathcal{R}_\mathrm{area}$   &  $\#$Fold-\\ 
& {(secs.)}  &                         &   SD     & ings                         & & {(secs.)}  &                        &   SD     & ings\\
\hline 
Right Hand   & 0.75  &$2.44 \times 10^{-2}$ &$6.77 \times 10^{-2}$ & 0  & & 0.81  &$8.55 \times 10^{-2}$ &$2.06 \times 10^{-1}$  &$165$\\
David Head   & 0.72  &$2.12 \times 10^{-3}$ &$1.28 \times 10^{-2}$ & 0  & & 0.93  &$3.65 \times 10^{-3}$ &$1.70 \times 10^{-2}$  &0\\
Bull         & 3.57  &$2.28 \times 10^{-2}$ &$5.24 \times 10^{-2}$ & 0  & & 3.73  &$2.05 \times 10^{-1}$ &$1.38 \times 10^{-1}$  &$ 18$\\
Bulldog      & 9.93  &$1.90 \times 10^{-3}$ &$1.57 \times 10^{-2}$ & 0  & & 4.59  &$3.01 \times 10^{-2}$ &$4.98 \times 10^{-2}$  &0\\
Lion Statue  & 15.73 &$5.81 \times 10^{-3}$ &$2.14 \times 10^{-2}$ & 0  & &14.60  &$9.33 \times 10^{-2}$ &$8.48 \times 10^{-2}$  &$ 11$\\
Gargoyle     & 15.97 &$3.83 \times 10^{-3}$ &$2.12 \times 10^{-2}$ & 0  & & 7.23  &$1.66 \times 10^{-2}$ &$3.89 \times 10^{-2}$  &0\\
Max Planck   & 4.73  &$9.51 \times 10^{-4}$ &$9.35 \times 10^{-3}$ & 0  & & 2.58  &$2.49 \times 10^{-2}$ &$4.48 \times 10^{-2}$  &0\\
Chess King   & 46.61 &$1.98 \times 10^{-2}$ &$3.90 \times 10^{-2}$ & 0  & &39.56  &$4.58 \times 10^{-2}$ &$7.27 \times 10^{-2}$  &$381$\\
\specialrule{.2em}{.1em}{.1em}
\multicolumn{9}{l}{\small $^\dagger$ Stopping criteria: energy deficit $< 10^{-5}$ or reach 100 iterations (Algorithm \ref{alg:PCG}). } \\
\end{tabular}
}
\end{table}

\begin{table}[ht]
\centering
\caption{  { Numerical results of Reimannian bijective correction (Algorithm~\ref{alg:bij_cor}) after minimizing $E_\rA$ \eqref{eq:Ea} in Table \ref{tab:Ea_comparison}.
$\#$Iter.: iteration count; $\mathcal{R}_\mathrm{area}$: area ratio \eqref{eq:Area_dist}; $\#$Foldings: number of folded triangles.}
}
\label{tab:bij_correct}
\resizebox{\textwidth}{!}{
\begin{tabular}{lrrrcc}
\specialrule{.2em}{.1em}{.1em}
\multirow{3}{*}{Model name} & \multicolumn{5}{c}{Riemannian bijective correction (Algorithm \ref{alg:bij_cor})}\\
\cline{2-6} 
& Time       & \multirow{2}{*}{$\#$Iter.} & $\#$Fold- & \multirow{2}{*}{$E_\rA$$^\dagger$ }  &$\mathcal{R}_\mathrm{area}$ \\ 
& {(secs.)}  &                             &  ings$^\dagger$     &                           & SD$^\dagger$  \\
\hline 
Right Hand   & 0.1   & 7 & $165 \rightarrow 0$ &$8.55 \times 10^{-2} \rightarrow 2.67 \times 10^{-1}$ & $9.57 \times 10^{-2} \rightarrow 2.06 \times 10^{-1}$ \\
Bull         & 0.14  & 7 & $18 \rightarrow 0$  &$2.05\times 10^{-1} \rightarrow 2.21 \times 10^{-1}$  & $1.27 \times 10^{-1} \rightarrow 1.38 \times 10^{-1}$ \\
Lion Statue  & 0.22  & 3 & $11 \rightarrow 0$  &$9.33 \times 10^{-2} \rightarrow 9.57 \times 10^{-2}$ & $8.37 \times 10^{-2} \rightarrow 8.48 \times 10^{-2}$ \\
Chess King   & 25.12 &47 & $381 \rightarrow 0$ &$4.58 \times 10^{-2} \rightarrow 6.89 \times 10^{-2}$ & $5.82 \times 10^{-2} \rightarrow 7.24 \times 10^{-2}$ \\
\specialrule{.2em}{.1em}{.1em}
\multicolumn{6}{l}{\small $^\dagger$ Left: without bijective correction; Right: with bijective correction.  } \\
\end{tabular}
}
\end{table}

\subsection{{ Comparison under different initial mappings}} \label{sec:6.2}
\begin{figure}[ht]
\centering
\resizebox{\textwidth}{!}{
\begin{tabular}{cc}
\includegraphics[height=4.5cm]{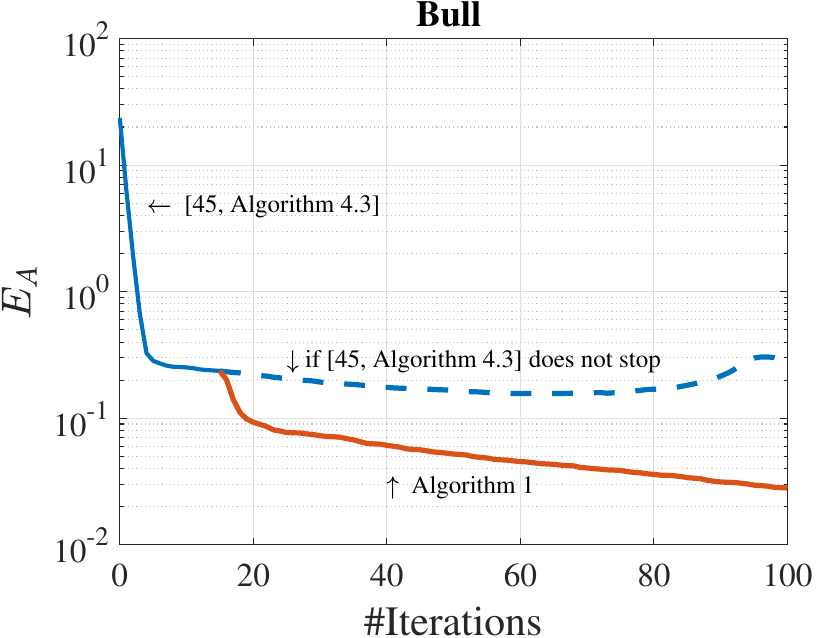} &
\includegraphics[height=4.5cm]{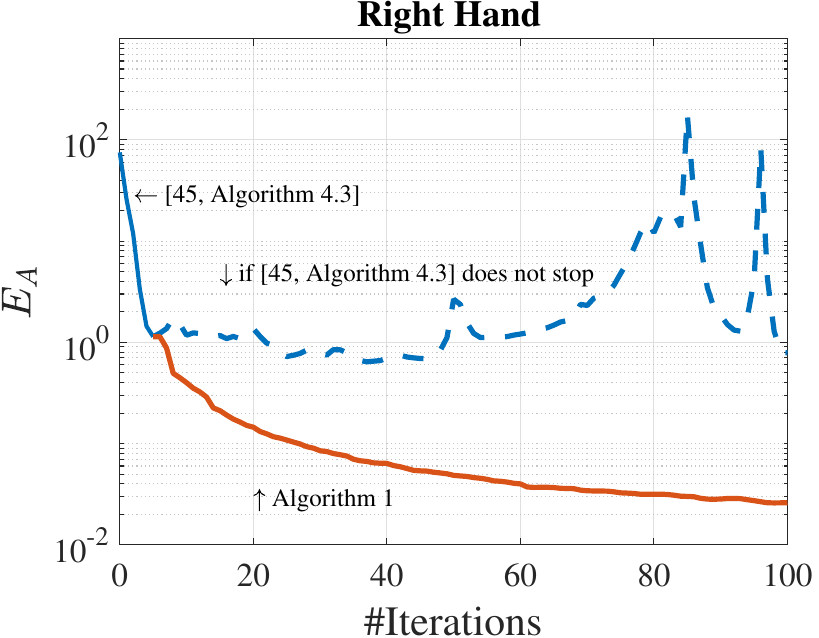} 
\end{tabular}
}
\caption{{ 
The illustration of the authalic energy $E_\rA$ versus the number of iterations. 
The solid blue line represents the energy $E_\rA$ of the initial mapping computed by the fixed-point method \cite[Algorithm 4.3]{YuLL19}. The solid orange line shows the energy decrease achieved by Algorithm~\ref{alg:PCG} over 100 iterations. For comparison, the dashed blue line illustrates the energy $E_\rA$ obtained by running the fixed-point method for 100 iterations without early termination.
}}
\label{fig:SEM_EnergyIter}
\end{figure}

\begin{figure}[ht]
\centering
\resizebox{\textwidth}{!}{
\begin{tabular}{cc}
\includegraphics[height=4.5cm]{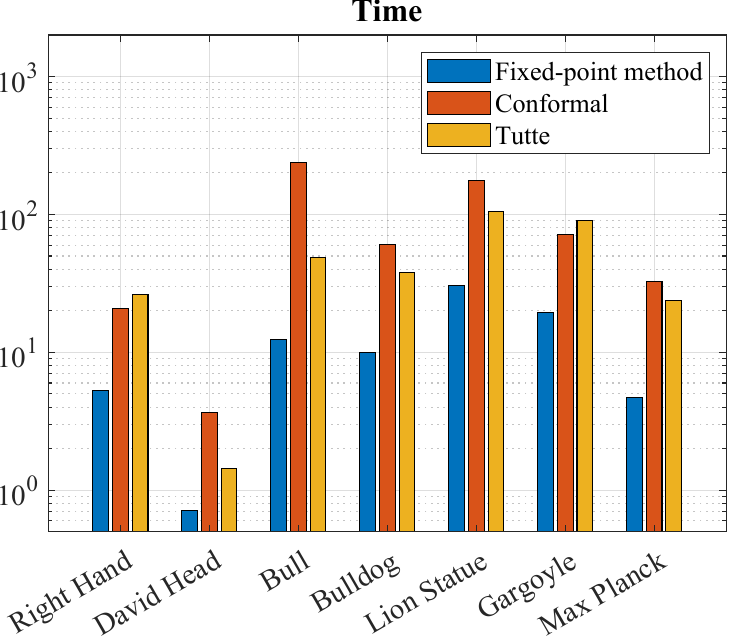} &
\includegraphics[height=4.5cm]{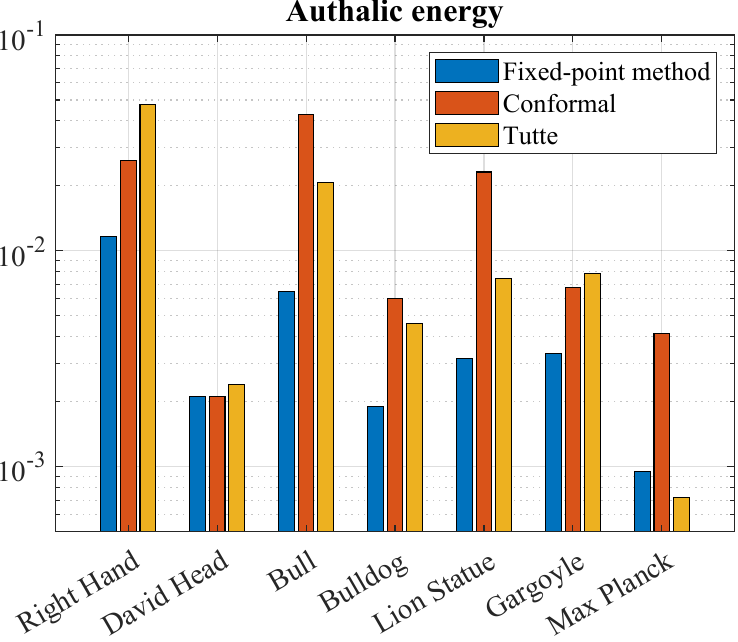} 
\end{tabular}
}
\caption{{ Numerical results of Algorithm \ref{alg:PCG} with different initialization: the fixed-point method \cite[Algorithm 4.3]{YuLL19}, the spherical conformal map \cite{HaAT00}, and the spherical Tutte map \cite{HaAT00}. Left: computational time cost; Right: Authalic energy $E_\rA$ in \eqref{eq:Ea}. }}
\label{fig:Initial}
\end{figure}

{

As mentioned in Section~\ref{sec:4.1}, we initialize the mapping using the fixed-point method \cite[Algorithm 4.3]{YuLL19}, running it for up to 15 iterations or terminating early if $E_\rA$ increases. Although this method provides a good initial mapping, it does not guarantee a decrease in $E_\rA$ at each step and may diverge if too many iterations are applied.

Figure~\ref{fig:SEM_EnergyIter} shows the authalic energy $E_\rA$ versus number of iterations for both the fixed-point method \cite[Algorithm 4.3]{YuLL19} and Algorithm \ref{alg:PCG}, applied to the Bull and Right Hand mesh models. The solid blue line represents the energy $E_\rA$ of the initial mapping computed by the fixed-point method. The solid orange line shows the energy decrease achieved by Algorithm~\ref{alg:PCG} over 100 iterations. For comparison, the dashed blue line illustrates the energy $E_\rA$ obtained by running the fixed-point method for 100 iterations without early termination.

We observe that the fixed-point method significantly reduces the energy $E_\rA$ in the first few iterations. However, its effectiveness diminishes rapidly and may even increase in energy when run for too many iterations. In contrast, the subsequent optimization using Algorithm~\ref{alg:PCG} consistently decreases the energy $E_\rA$ throughout the 100 iterations, which is consistent with the convergence guarantee.

To evaluate the importance of the initial mapping, we perform the numerical experiments of Algorithm~\ref{alg:PCG} (terminated when the energy deficit falls below $10^{-5}$) with three different initial mappings: the fixed-point method \cite[Algorithm 4.3]{YuLL19}, the spherical conformal map \cite{HaAT00}, and the spherical Tutte map \cite{HaAT00} (computed by replacing the cotangent Laplacian with the Tutte Laplacian). The Chess King model is excluded due to inefficiencies encountered when using the conformal and Tutte initial maps.

The numerical results are summarized in Table~\ref{tab:Initial_selection} and visualized in Figure~\ref{fig:Initial}. As expected, the mappings resulting from the fixed-point initial map achieve the highest accuracy and efficiency across most mesh models. While the other initial maps still lead to reasonably good final maps, as indicated by the values of $E_\rA$, they generally require significantly larger computational time costs, such as the Bull model when using the conformal initial map. We also observe that poorer initial maps tend to produce more folded triangles in the resulting mappings. Nevertheless, these folded triangles are successfully eliminated by applying the Riemannian bijective correction.

In summary, although Algorithm \ref{alg:PCG} demonstrates robustness to different initial mappings due to its global convergence property, the fixed-point method \cite[Algorithm 4.3]{YuLL19} provides a good initial map for achieving both high accuracy and computational efficiency.
}

\begin{table}[ht]
\centering
\caption{ { Numerical results of Algorithm~\ref{alg:PCG} (terminated when the energy deficit falls below $10^{-5}$) with three different initial mappings: the fixed-point method \cite[Algorithm 4.3]{YuLL19}, the spherical conformal map \cite{HaAT00}, and the spherical Tutte map \cite{HaAT00}.
$E_\rA$: authalic energy \eqref{eq:Ea}; $\#$Foldings: number of folded triangles. }
}
\label{tab:Initial_selection}
\resizebox{\textwidth}{!}{
\begin{tabular}{lrcrcrcrcrcr}
\specialrule{.2em}{.1em}{.1em}
\multirow{3}{*}{Model name} & \multicolumn{3}{c}{Fixed-point method \cite{YuLL19} } & & \multicolumn{3}{c}{Conformal \cite{HaAT00} $^\dagger$} & & \multicolumn{3}{c}{Tutte \cite{HaAT00}$^\ddagger$ } \\
\cline{2-4} \cline{6-8} \cline{10-12}
& Time   & \multirow{2}{*}{$E_\rA$} & $\#$Fold- & & Time  & \multirow{2}{*}{$E_\rA$} & $\#$Fold- & & Time   & \multirow{2}{*}{$E_\rA$} & $\#$Fold- \\
& {(secs.)}  &    & ings  & & {(secs.)}  &    & ings$^*$ & & {(secs.)}  &    & ings$^*$ \\
\hline 
Right Hand   & 5.30  &$1.17 \times 10^{-2}$ & 0 & & 20.70 &$2.61 \times 10^{-2}$ & $3 \rightarrow0$ & & 26.29 &$4.77 \times 10^{-2}$ & $11\rightarrow0$\\
David Head   & 0.72  &$2.12 \times 10^{-3}$ & 0 & & 3.65  &$2.12 \times 10^{-3}$ & 0                & & 1.44  &$2.40 \times 10^{-3}$ & 0\\
Bull         & 12.48 &$6.48 \times 10^{-3}$ & 0 & & 236.25 &$4.29 \times 10^{-2}$ & $7\rightarrow0$ & & 48.73 &$2.07 \times 10^{-2}$ & $5 \rightarrow0$\\
Bulldog      & 9.93  &$1.90 \times 10^{-3}$ & 0 & & 60.94 &$6.02 \times 10^{-3}$ & $7 \rightarrow0$ & & 38.18 &$4.59 \times 10^{-3}$ & $5 \rightarrow0$\\
Lion Statue  & 30.75 &$3.18 \times 10^{-3}$ & 0 & & 177.04&$2.32 \times 10^{-2}$ & $6 \rightarrow0$ & & 105.51&$7.48 \times 10^{-3}$ & $7 \rightarrow0$\\
Gargoyle     & 19.60 &$3.35 \times 10^{-3}$ & 0 & & 71.38 &$6.80 \times 10^{-3}$ & $7 \rightarrow0$ & & 90.01 &$7.84 \times 10^{-3}$ & 0\\
Max Planck   & 4.73  &$9.51 \times 10^{-4}$ & 0 & & 32.87 &$4.14 \times 10^{-3}$ & $5 \rightarrow0$ & & 23.72 &$7.19 \times 10^{-4}$ & 0\\
\specialrule{.2em}{.1em}{.1em}
\multicolumn{12}{l}{\small $^*$ Left: without bijective correction; Right: with bijective correction.} \\
\multicolumn{12}{l}{\small $^\dagger$ The subsequent Algorithm~\ref{alg:PCG} uses the cotangent-weighted Laplacian matrix on the image as the preconditioner.} \\
\multicolumn{12}{l}{\small $^\ddagger$ The Tutte Laplacian is used in constructing both the map and the preconditioner in Algorithm~\ref{alg:PCG}.} 
\end{tabular}
}
\end{table}

\subsection{Comparison with state-of-the-art methods} \label{sec:6.3}
Next, we compare our proposed method (SAEM) with the spherical density-equalizing mapping (SDEM) \cite{LyZh24} and Riemannian gradient descent (RGD) method \cite{SuYu24}. As reported in \cite{SuYu24}, RGD achieves superior area preservation compared to the fixed-point method \cite{YuLL19}, adaptive area-preserving mapping \cite{ChGK22}, and optimal transportation mapping \cite{CuQW19}.

\begin{figure}[ht]
\centering
\resizebox{\textwidth}{!}{
\begin{tabular}{cc}
\includegraphics[height=4.5cm]{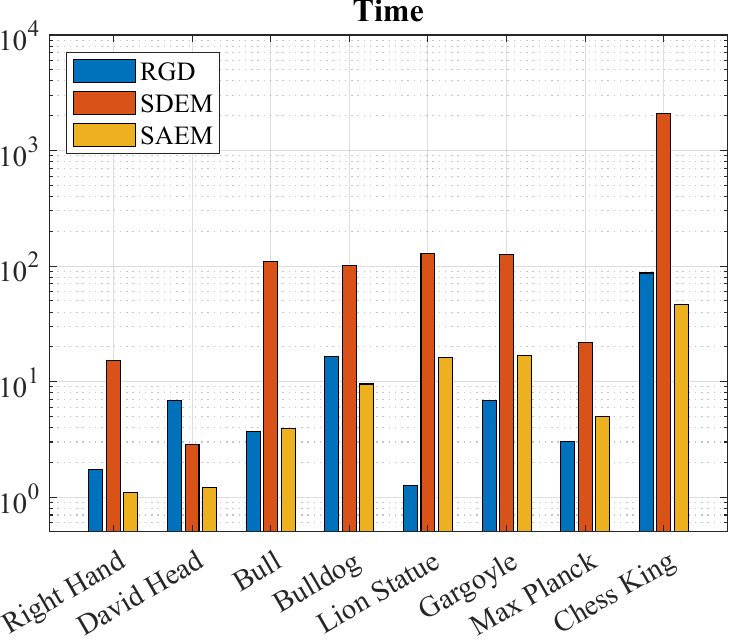} &
\includegraphics[height=4.5cm]{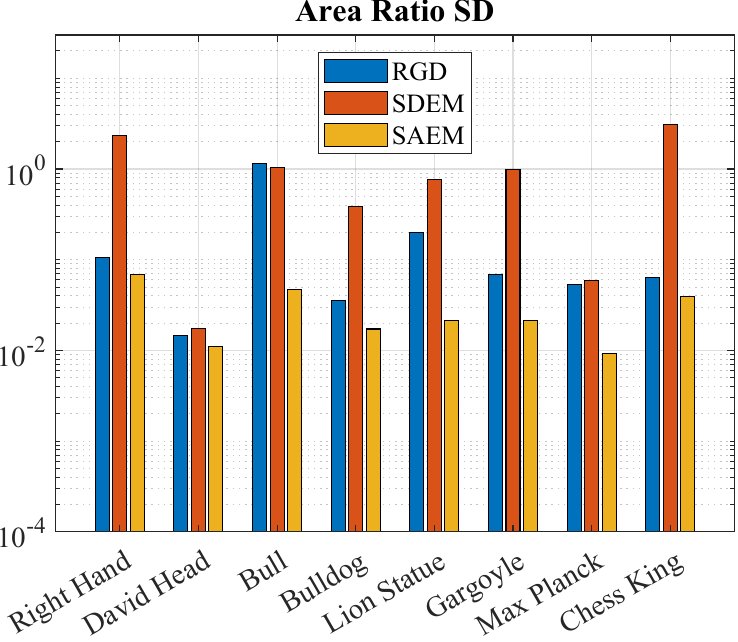} 
\end{tabular}
}
\caption{Numerical results of our proposed method (SAEM), Riemannian gradient descent (RGD) \cite{SuYu24}, and spherical density-equalizing mapping (SDEM) \cite{LyZh24}. Left: computational time cost; Right: standard deviation of area ratio \eqref{eq:Area_dist}.}
\label{fig:SD_Time}
\end{figure}

The numerical results are summarized in Table~\ref{tab:method_comparison}. In terms of computational efficiency, our proposed method outperforms the compared methods in most cases (see Figure~\ref{fig:SD_Time}, left). The inefficiency of SDEM is likely caused by solving a large-scale linear system for the PDE and bijective correction in each iteration. Meanwhile, both our method and RGD use line searches, leading to comparable efficiency.
Regarding accuracy, both $E_\rA$ and the standard deviation of $\mathcal{R}_\mathrm{area}$ show that our method achieves superior area-preserving performance (see Figure~\ref{fig:SD_Time}, right). The relatively poor accuracy of the RGD method may be attributed to the intrinsic limitations of gradient descent. 

Furthermore, our Riemannian bijective correction yields zero folded triangles in all SAEM maps. The manipulation of the Beltrami coefficient in SDEM eliminates foldings in most cases, but leaves four folded triangles for Lion and one for Gargoyle. In contrast, the bijective correction step in the RGD method fails to eliminate all folded triangles in certain challenging cases, possibly due to truncation errors by the stereographic projection.

\begin{figure}[ht]
\centering
\includegraphics[width=0.8\textwidth]{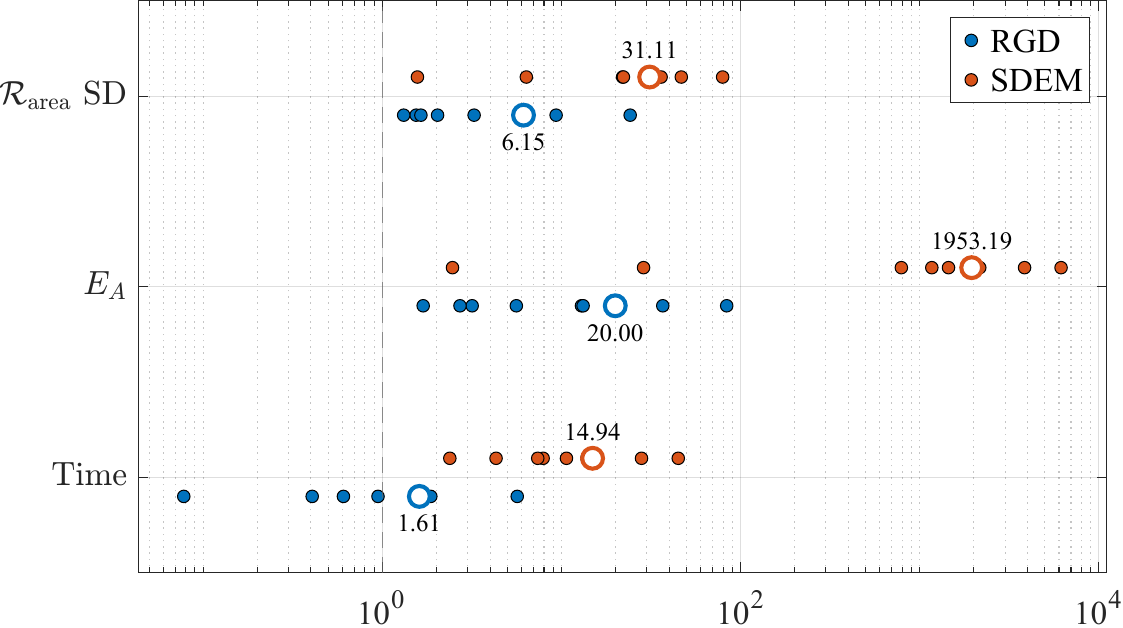}
\caption{
The ratio of the results from Riemannian gradient descent (RGD) \cite{SuYu24} and spherical density-equalizing mapping (SDEM) \cite{LyZh24} to those of our method (SAEM). The big circle indicates the average of the ratios.
$E_\rA$: authalic energy \eqref{eq:Ea}; $\mathcal{R}_\mathrm{area}$: area ratio \eqref{eq:Area_dist}.}
\label{fig:method_Ranking}
\end{figure}

To clearly demonstrate the improvements of our proposed method, Figure~\ref{fig:method_Ranking} presents the ratios of the results from RGD and SDEM to those of our method across the metrics in Table \ref{tab:method_comparison}, including computational time, authalic energy $E_\rA$, and the standard deviation of $\mathcal{R}_\mathrm{area}$. A ratio of 1 indicates equivalent performance, while a ratio greater than 1 signifies the superiority of our algorithm. The results emphasize the outstanding performance of our method compared to state-of-the-art methods.

In summary, our proposed method outperforms the RGD method and SDEM in terms of efficiency, accuracy, and bijectivity.

\begin{table}[p]
\centering
\caption{Numerical results of spherical authalic energy minimization (SAEM), Riemannian gradient descent (RGD) \cite{SuYu24}, and spherical density-equalizing mapping (SDEM) \cite{LyZh24}.
$E_\rA$: authalic energy \eqref{eq:Ea}; $\mathcal{R}_\mathrm{area}$: area ratio \eqref{eq:Area_dist}; $\#$Foldings: number of folding triangles; $\#$Iterations: iteration count.
}
\label{tab:method_comparison}
\resizebox{\textwidth}{!}{
\begin{tabular}{clccc}
\specialrule{.2em}{.1em}{.1em}
\multirow{2}{*}{Metric} & \multirow{2}{*}{Model name}  & Spherical authalic energy & Riemannian gradient  & Spherical density-equalizing\\
 & & minimization (SAEM)$^\dagger$  & descent (RGD)$^\ddagger$ \cite{SuYu24}  &  mapping (SDEM)$^\mathsection$ \cite{LyZh24}\\
\hline
\multirow{8}{*}{Time} 
& Right Hand         & 1.09       & 1.73        & 15.20\\
& David Head         & 1.20       & 6.81        &  2.86\\
& Bull               & 3.91       & 3.70        &109.48\\
& Bulldog            & 9.53       & 16.55       &101.73\\
& Lion               &16.16       & 1.26        &127.85\\
& Gargoyle           &16.93       & 6.88        &124.97\\
& Max Planck         & 5.02       & 3.05        & 21.70\\
& Chess King         &46.61       &87.03        &2094.42\\
\hline
\multirow{8}{*}{$E_\rA$} 
& Right Hand         &$2.56 \times 10^{-2}$  &$8.12 \times 10^{-2}$  &$5.54 \times 10^{1}$\\
& David Head         &$1.54 \times 10^{-3}$  &$2.60 \times 10^{-3}$  &$3.79 \times 10^{-3}$\\
& Bull               &$1.85 \times 10^{-2}$  &$2.40 \times 10^{-1}$  &$2.17 \times 10^{1}$\\
& Bulldog            &$2.34 \times 10^{-3}$  &$1.31 \times 10^{-2}$  &$1.85 \times 10^{0}$\\
& Lion               &$5.72 \times 10^{-3}$  &$4.79 \times 10^{-1}$  &$8.29 \times 10^{0}$\\
& Gargoyle           &$3.89 \times 10^{-3}$  &$5.14 \times 10^{-2}$  &$1.50 \times 10^{1}$\\
& Max Planck         &$9.61 \times 10^{-4}$  &$3.54 \times 10^{-2}$  &$2.76 \times 10^{-2}$\\
& Chess King         &$1.98 \times 10^{-2}$  &$5.38 \times 10^{-2}$  &$1.22 \times 10^{2}$\\
\hline
\multirow{8}{*}{$\mathcal{R}_\mathrm{area}$ SD} 
& Right Hand         &$6.79 \times 10^{-2}$  &$1.05 \times 10^{-1}$  &$2.35 \times 10^{0}$\\
& David Head         &$1.10 \times 10^{-2}$  &$1.45 \times 10^{-2}$  &$1.73 \times 10^{-2}$\\
& Bull               &$4.71 \times 10^{-2}$  &$1.14 \times 10^{~0}$  &$1.03 \times 10^{0}$\\
& Bulldog            &$1.72 \times 10^{-2}$  &$3.51 \times 10^{-2}$  &$3.83 \times 10^{-1}$\\
& Lion               &$2.14 \times 10^{-2}$  &$2.00 \times 10^{-1}$  &$7.71 \times 10^{-1}$\\
& Gargoyle           &$2.11 \times 10^{-2}$  &$6.88 \times 10^{-2}$  &$9.87 \times 10^{-1}$\\
& Max Planck         &$9.26 \times 10^{-3}$  &$5.38 \times 10^{-2}$  &$5.90 \times 10^{-2}$\\
& Chess King         &$3.90 \times 10^{-2}$  &$6.41 \times 10^{-2}$  &$3.10 \times 10^{0}$\\
\hline
\multirow{8}{*}{$\#$ Foldings} 
& Right Hand         & 0       & 0      & 0\\
& David Head         & 0       & 0      & 0\\
& Bull               & 0       & 3      & 0\\
& Bulldog            & 0       & 0      & 0\\
& Lion               & 0       & 0      & 4\\
& Gargoyle           & 0       & 0      & 1\\
& Max Planck         & 0       & 0      & 0\\
& Chess King         & 0       &17      & 0\\
\hline
\multirow{8}{*}{$\#$ Iterations} 
& Right Hand         & $100^*$ & 127    & $50^*$\\
& David Head         & 26       &$200^*$ & 26 \\
& Bull               & $100^*$ & 50     & $50^*$\\
& Bulldog            & 50      & 60     & $50^*$\\
& Lion               & $100^*$ & 2      & $50^*$\\
& Gargoyle           & $100^*$ & 17     & $50^*$\\
& Max Planck         & 24      & 8      & 23 \\
& Chess King         & $100^*$ & 107    & $50^*$\\
\specialrule{.2em}{.1em}{.1em}
\multicolumn{5}{l}{\small $^*$ indicates reaching the maximum iteration limit.} \\
\multicolumn{5}{l}{\small $^\dagger$ stopping criteria: energy deficit $< 10^{-5}$ or reach 100 iterations.} \\
\multicolumn{5}{l}{\small $^\ddagger$ stopping criteria: energy deficit $< 5 \times 10^{-6}$ or reach 200 iterations.} \\
\multicolumn{5}{l}{\small $^\mathsection$ stopping parameter $= 10^{-3}$, maximum iterations $= 50$.} \\
\end{tabular}
}
\end{table}

\section{Application in shape description}
\label{sec:7}
{ \em Spherical harmonics} \cite{Byer59} are a set of orthogonal basis functions defined on the unit sphere $\mathbb{S}^2$. { Any function on the unit sphere can be represented as a linear combination of spherical harmonics.} 
Consequently, spherical harmonics can be used to effectively represent surfaces through spherical parameterization \cite{BrGK95}.

{
In particular, once a spherical parameterization $f:\M\to\mathbb{S}^2$ is established, the surface $\M$ can be viewed as the image of the inverse mapping $f^{-1}$ defined on the unit sphere $\mathbb{S}^2$. That is, each vertex $v\in\V(\M)$ is represented as $v(\theta, \phi)$.} 

Then, the vertex $v$ can be represented as a linear combination of spherical harmonic basis functions $Y_\ell^m (\theta, \phi)$ of degree $\ell$ and order $m \in [-\ell, \ell]$ as
\[
v(\theta, \phi) = \sum_{\ell = 0}^\infty \sum_{m = -\ell}^\ell \mathbf{c}_\ell^m Y_\ell^m(\theta, \phi),
\]
where $\mathbf{c}_\ell^m$ is the $3 \times 1$ vector of complex coefficients corresponding to each coordinate of $v_i$.
The spherical harmonic $Y_\ell^m$ is a complex-valued function defined by 
\begin{equation*}
Y_\ell^m(\theta, \phi) = N_\ell^m  e^{i m \phi} P_\ell^m (\cos \theta),
\end{equation*}
where $P_\ell^m$ are the { \textit{associated Legendre polynomials}},
\[
P_\ell^m(x) = \frac{(-1)^m}{2^\ell \,\ell!} (1-x^2)^{m/2} \frac{\d^{m+\ell}}{\d x^{m+\ell}} (x^2 - 1)^\ell,
\]
and $N_\ell^m$ is a normalization constant given by
\[
N_\ell^m = \sqrt{\frac{2\ell + 1}{4 \pi} \frac{(\ell-m)!}{(\ell+m)!}}.
\]

An important property of these polynomials is the orthogonality relation
\[
\int_{-1}^1 P_\ell^m(x) P_{\ell'}^m(x) \,\d x = \frac{2(m+\ell)!}{(2\ell + 1)(\ell-m)!} \delta_{\ell, \ell'},
\]
where $x = \cos \theta \in [-1,1]$, and $\delta$ denotes the { \textit{Kronecker delta function}, given by
\[
\delta_{i,j} = 
\begin{cases}
    \,1 & \mbox{if}~ i = j,\\
    \,0 & \mbox{if}~ i \neq j.
\end{cases}
\]

}

\begin{figure}[tbp]
\centering
\resizebox{\textwidth}{!}{
\begin{tabular}{cccc}
$d=1$ & $d=4$ & $d=9$ & $d=16$ \\
\includegraphics[height=4cm]{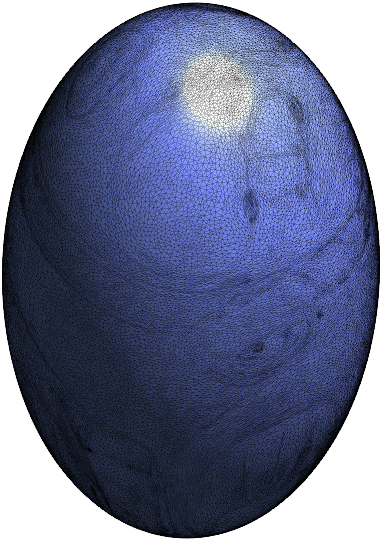} &
\includegraphics[height=4cm]{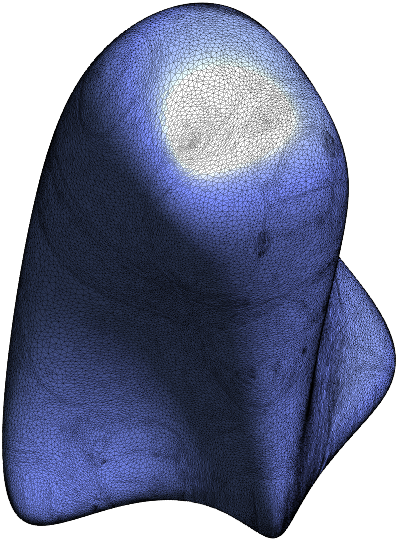} &
\includegraphics[height=4cm]{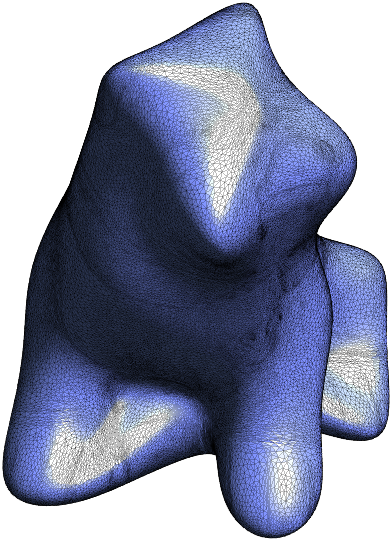} &
\includegraphics[height=4cm]{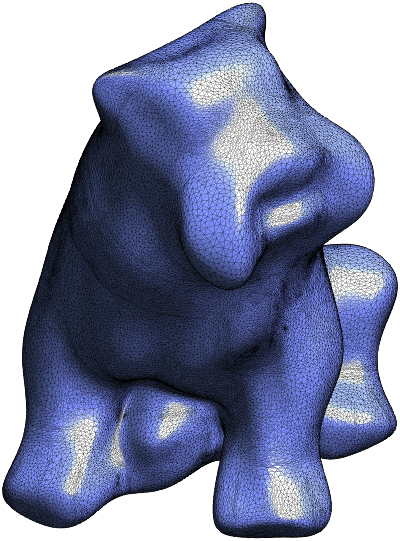} \\ 
$d=25$ & $d=36$ & $d=49$ & $d=64$ \\
\includegraphics[height=4cm]{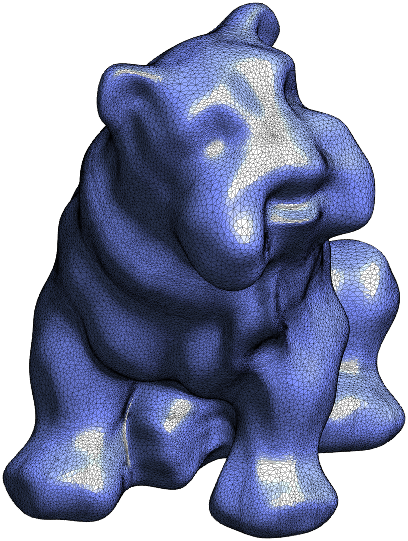} &
\includegraphics[height=4cm]{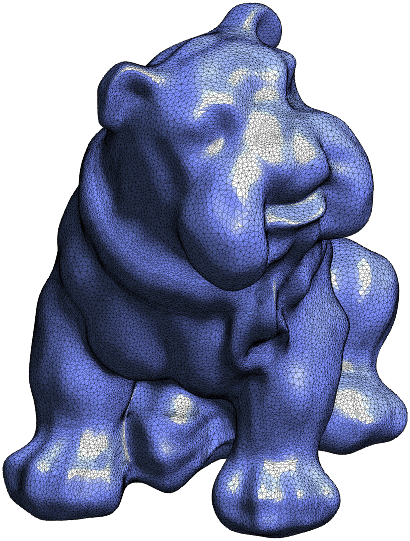} &
\includegraphics[height=4cm]{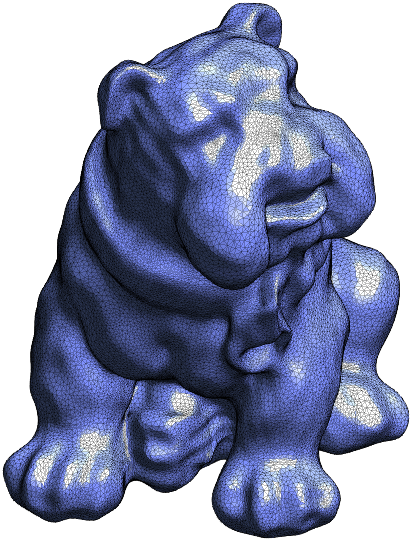} &
\includegraphics[height=4cm]{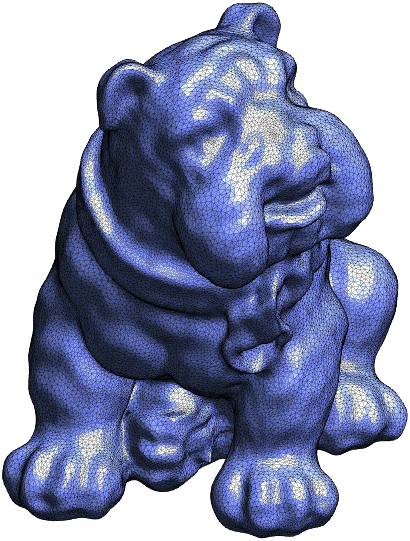}  
\end{tabular}
}
\caption{Approximation of the Bulldog triangular mesh using spherical harmonics in \eqref{eq:SH_approx_matrix} with degrees $d = 1, 4, 9, 16, 25, 36, 49,$ and $64$, respectively.}
\label{fig:SH_Bulldog}
\end{figure}

In practice, we truncate the infinite summation at a finite degree $d$, thereby approximating each vertex $v_i$ as 
\begin{equation}
v(\theta, \phi) \approx \sum_{\ell = 0}^d \sum_{m = -\ell}^\ell \mathbf{c}_\ell^m Y_\ell^m(\theta, \phi).
\label{eq:SH_approx}
\end{equation}
Let $\mathbf{V} = (v_1, v_2, \ldots, v_n)^\top$ be $n \times 3$ matrix of vertex coordinates, and let $\bC$ be the $(d+1)^2 \times 3$ matrix of coefficients. Define $\bY$ as the $n \times (d+1)^2$ matrix by
\[
\bY = 
\begin{bmatrix}
    Y_0^0(\theta_1, \phi_1)  &  Y_1^{-1}(\theta_1, \phi_1)  &  Y_1^{0}(\theta_1, \phi_1)  &  Y_1^{1}(\theta_1, \phi_1)  &  \cdots  & Y_d^{d}(\theta_1, \phi_1)\\
    \vdots  &  \vdots  &  \vdots  &  \vdots  &  \vdots  &  \vdots\\
    Y_0^0(\theta_n, \phi_n)  &  Y_1^{-1}(\theta_n, \phi_n)  &  Y_1^{0}(\theta_n, \phi_n)  &  Y_1^{1}(\theta_n, \phi_n)  &  \cdots  & Y_d^{d}(\theta_n, \phi_n)
\end{bmatrix}.
\]
Then, \eqref{eq:SH_approx} can be expressed in matrix form as
\begin{equation}
\mathbf{V} \approx \bY \, \bC.
\label{eq:SH_approx_matrix}
\end{equation}
The coefficient matrix $\bC$ can be estimated by the least-square method as 
\[
\bC = (\bY^\top \bY)^{-1} \bY^\top \mathbf{V}.
\]
Therefore, once the spherical parameterization is computed by Algorithm~\ref{alg:PCG}, the surface can be represented by the coefficient matrix $\bC$ up to degree $d$.
Figure~\ref{fig:SH_Bulldog} presents the surface approximation of the Bulldog triangular mesh as in \eqref{eq:SH_approx_matrix}, with degrees 1, 4, 9, 16, 25, 36, 49, and 64. As expected, increasing the degree of the spherical harmonics yields progressively finer and more detailed reconstructions.

\begin{figure}[tbp]
\centering
\resizebox{\textwidth}{!}{
\begin{tabular}{cccc}
\multicolumn{2}{c}{Conformal} & \multicolumn{2}{c}{Authalic} \\
\includegraphics[height=4cm]{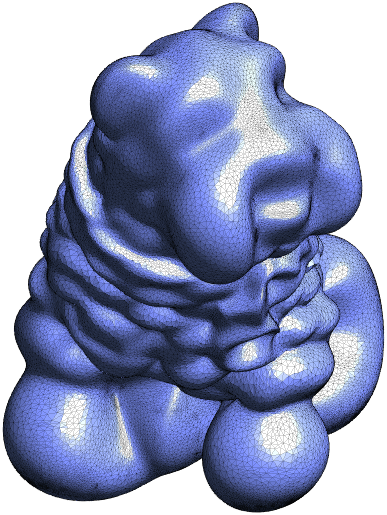} &
\includegraphics[height=4cm]{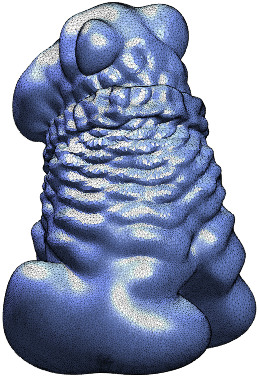} &
\includegraphics[height=4cm]{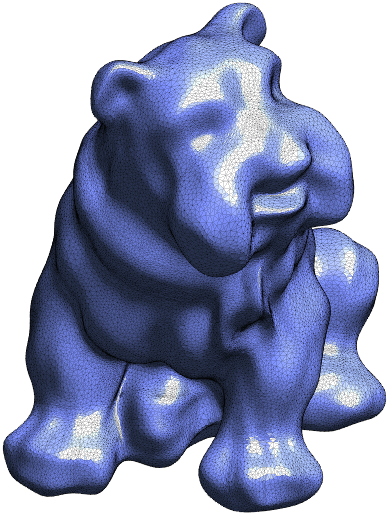} &
\includegraphics[height=4cm]{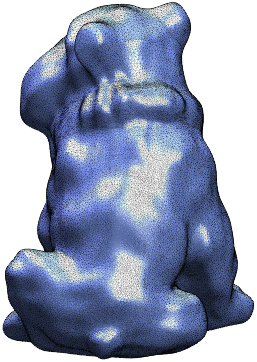} 
\end{tabular}
}
\caption{Comparison of conformal and authalic parameterizations for reconstructing the Bulldog triangular mesh using spherical harmonics in \eqref{eq:SH_approx_matrix} with degree 30.}
\label{fig:SH_Bulldog_CEM}
\end{figure}

\begin{figure}[tbp]
\centering
\resizebox{\textwidth}{!}{
\begin{tabular}{cccc}
\multicolumn{2}{c}{Conformal} & \multicolumn{2}{c}{Authalic} \\
\includegraphics[height=4cm]{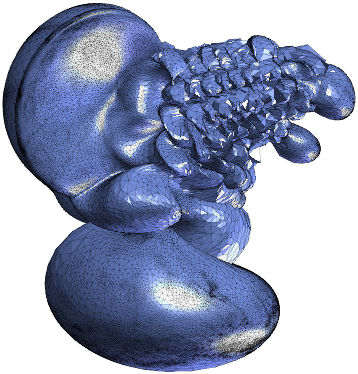} &
\includegraphics[height=4cm]{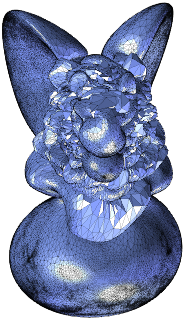} &
\includegraphics[height=4cm]{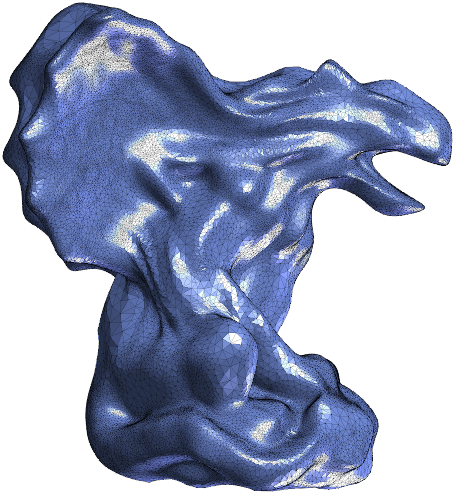} &
\includegraphics[height=4cm]{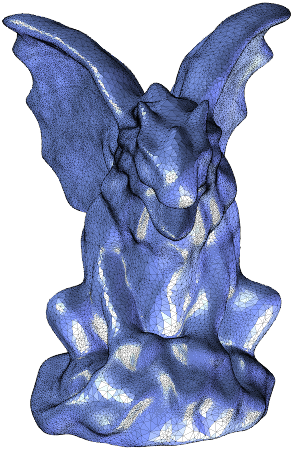} 
\end{tabular}
}
\caption{Comparison of conformal and authalic parameterizations for reconstructing the Gargoyle triangular mesh using spherical harmonics in \eqref{eq:SH_approx_matrix} with degree 30.}
\label{fig:SH_Gargoyle_CEM}
\end{figure}

Area preservation in spherical parameterization is essential for the effectiveness of spherical harmonics in shape description. We compare surface approximation results on the Bulldog mesh using authalic parameterization (Algorithm~\ref{alg:PCG}) and conformal parameterization \cite{ChLL15} as illustrated in Figure~\ref{fig:SH_Bulldog_CEM}. The authalic parameterization yields accurate surface reconstruction, while the conformal parameterization performs poorly due to its lack of uniform vertex distribution on the unit sphere.
The limitations of conformal parameterization become even more pronounced for geometrically complicated surfaces, as evident in the Gargoyle mesh shown in Figure~\ref{fig:SH_Gargoyle_CEM}.

\section{Discussion}
\label{sec:8}

In this section, we highlight the key differences between the proposed method and existing approaches.

\subsection{Computational methods for spherical authalic mappings}
There are several strategies for computing spherical area-preserving mappings. One common approach is to treat the sphere as a plane via stereographic projection. For example, Nadeem et al.~\cite{NaSZ17} and Choi et al.~\cite{ChGK22} obtain the area-preserving mapping through optimal mass transport theory with the stereographic projection. They start from an initial spherical conformal mapping and then solve the problem using Newton's method on the stereographic plane. Similarly, Yueh et al.~\cite{YuLL19} employ stereographic projection and then apply the fixed-point method to update the unit disk for the southern hemisphere while using inversion for the northern hemisphere to obtain spherical area-preserving mappings.

A second approach is based on Riemannian optimization. Variables are updated on the local coordinates by projecting onto the tangent plane, and then retracted to the unit sphere, in line with the Riemannian geometry as a generalization of Euclidean space. Lyu et al.~\cite{LyZh24} use the area ratio as a density function and achieve area preservation by solving the diffusion equation on the local tangent plane. Sutti and Yueh~\cite{SuYu24} apply the Riemannian gradient descent method for energy minimization, which offers better theoretical convergence properties than the projected gradient descent method.

The final approach addresses this issue by using the global coordinate of the unit sphere. Cui et al.~\cite{CuQW19} compute the optimal transportation map using Newton's method in the spherical power diagram. In our work, we implement a preconditioned nonlinear CG method in spherical coordinates. 

\subsection{Spherical bijective correction}
{
In the work of Lyu et al.~\cite{LyZh24}, the Beltrami coefficients $\mu$ of a spherical map are computed on the stereographic plane. The coefficient $\mu$ is then manipulated to unfold triangles, relying on the fact that $\|\mu\|_\infty < 1$ ensures bijectivity for continuously differentiable mappings. However, since $\mu$ is updated at every iteration, the method requires solving a large-scale linear system repeatedly, which can be computationally inefficient.
}

In planar regions with convex boundaries, the unfolding post-processing can be performed by solving a linear system \cite{Yueh23} using a mean value weighted Laplacian matrix \cite{Floa03b}, as the convex combination map ensures bijectivity. To extend this planner method to spherical mappings, Sutti and Yueh \cite{SuYu24} combine it with the stereographic projection. While this approach successfully resolves folded triangles in the plane, it often fails to unfold them on the sphere due to truncation errors introduced by the stereographic projection. In contrast, Yueh et al. \cite{YuLL19} address this issue by locally resolving each folded triangle on the sphere. Because convex combination mappings are bijective in planar regions, we modify this method by resolving folded triangles on the tangent plane instead.

\subsection{SEM-based methods}
The stretch energy minimization (SEM) for disk-shaped area-preserving mappings was first introduced by Yueh et al. \cite{YuLW19} using the fixed-point method. This approach was later extended to spherical mappings via stereographic projection \cite{YuLL19}. However, from a numerical perspective, this method lacks global convergence, and in practice, the energy may increase for some challenging cases, as shown in Section \ref{sec:6.2}.

Inspired by this, Sutti and Yueh \cite{SuYu24} used the fixed-point method \cite{YuLL19} as an initial mapping, followed by the Riemannian gradient descent method to improve convergence. Although this modification resolved convergence issues, gradient descent remains intrinsically inefficient. Moreover, the objective functional that includes the inverse image area term may lead to triangle folding during optimization (see Section \ref{sec:3.1} for details).

In this work, we approximate the image area by the associated volume, thereby significantly improving the bijectivity of the resulting mappings. Furthermore, our preconditioned nonlinear CG method offers greater efficiency compared to the Riemannian gradient descent method.

\subsection{Parameterization in shape description}
Several studies have shown that area preservation of a spherical parameterization is essential for shape description with spherical harmonics, as first reported by  Brechb\"uhler et al. \cite{BrGK95}.
Similarly, Giri et al. \cite{GiCK21} demonstrated that area-preserving hemispherical parameterizations yield more accurate reconstructions compared to conformal spherical parameterizations. Likewise, Choi and Shaqfa \cite{ChSh25} further showed that area-preserving methods outperform conformal approaches in hemispheroidal parameterizations for shape reconstruction. 
{ Moreover, the quasi-conformal parameterization in \cite{ChLS20} is employed for shape description for early Alzheimer's disease detection with machine learning.}

Our findings are consistent with these prior studies. Authalic parameterizations produce better reconstruction results than conformal parameterizations, particularly on surfaces with complicated shapes. To the best of our knowledge, our proposed method achieves the best area preservation among state-of-the-art methods, indicating that it could be an ideal candidate for shape description via spherical harmonics.

\section{Conclusion}
\label{sec:9}

In this paper, we propose a novel objective functional, termed the spherical authalic energy, derived from the authalic energy and designed to enhance bijectivity during the minimization process for spherical area-preserving parameterizations. To effectively find a minimizer, we develop a preconditioned nonlinear CG method with theoretically guaranteed convergence. Numerical experiments show that our algorithm significantly improves both efficiency and accuracy compared to state-of-the-art methods. To further ensure bijectivity of the resulting mappings, we propose a Riemannian bijective correction as an optional post-processing step, supported by theoretical guarantees under mild assumptions. Finally, we demonstrate the practical utility of our method in shape description via spherical harmonics.

\subsection*{Acknowledgments}
The work of the authors was partially supported by the National Science and Technology Council and the National Center for Theoretical Sciences, Taiwan.

\subsection*{Data availability}
The MATLAB demo can be found at \url{https://web.ntnu.edu.tw/~81240002s/program/SAEM.zip}.

\footnotesize

\bibliographystyle{abbrv}

\bibliography{reference}

\end{document}